\documentclass{amsart}
\usepackage{amsfonts,amssymb,amscd,amsmath,enumerate,verbatim}

\usepackage[all]{xy}
\SelectTips{eu}{}
\xyoption{curve}

\newcommand{\xla}{\xleftarrow}
\newcommand{\lra}{\longrightarrow}

\newcommand{\xra}{\xrightarrow}

\newcommand{\card}{{\operatorname{card}}}

\newcommand{\op}{{{}^{\mathsf o}}}

\newcommand{\les}{{\scriptscriptstyle\leqslant}}

\newcommand{\shift}{{\sf\Sigma}}

\newcommand\col{\colon}
\newcommand{\cone}[1]{\mathsf{cone}(#1)}
\newcommand\dd{\partial}
\newcommand{\hh}[1]{{\operatorname{H}(#1)}}
\newcommand{\hH}{{\operatorname{H}}}
\newcommand{\HH}[2]{{\operatorname{H}_{#1}(#2)}}
\newcommand\BB{\operatorname{B}}

\newcommand\ZZ{\operatorname{Z}}

\newcommand{\Hom}[3]{\operatorname{Hom}_{#1}(#2,#3)}
\newcommand{\End}[2]{{\operatorname{End}_{#1}({#2})}}

\newcommand\Ker{\operatorname{Ker}}
\newcommand\Coker{\operatorname{Coker}}
\newcommand\Image{\operatorname{Im}}

\newcommand\Supp{\operatorname{Supp}}
\newcommand\gr[3]{{\operatorname{E}^{#1}\{{#2}^{#3}\}}}

\newcommand\Gr[5]{{{}^{#2}\!\operatorname{E}^{#1}_{#3}\{{#4}{}^{#5}\}}}
\newcommand\eGr[3]{{{}^{#2}\!\operatorname{E}^{#1}_{#3}}}

\newcommand\height{\operatorname{height}}
\newcommand\colim{\operatorname{colim}}

\newcommand\supheight{\operatorname{super\,height}}
\newcommand\length{\operatorname{length}}

\newcommand{\rank}[1]{\operatorname{rank}_{#1}}
\newcommand{\detrank}[1]{\operatorname{det\,rank}_{#1}}
\newcommand{\innrank}[1]{\operatorname{inn\,rank}_{#1}}

\newcommand{\Ann}{\operatorname{Ann}}
\newcommand\idmap{\operatorname{id}}

\newcommand\pd{\operatorname{proj\,dim}}

\newcommand\fm{{\mathfrak m}}

\newcommand\fp{{\mathfrak p}}
\newcommand\fq{{\mathfrak q}}

\newcommand\bsx{{\boldsymbol x}}

\newcommand\BF{{\mathbb F}}
\newcommand\BN{{\mathbb N}}
\newcommand\BQ{{\mathbb Q}}
\newcommand\BR{{\mathbb R}}
\newcommand\BZ{{\mathbb Z}}

\newcommand{\eps}{{\varepsilon}}

\newcommand{\var}{{\hskip1pt\vert\hskip1pt}}

\theoremstyle{plain}

\newtheorem{theorem}{Theorem}[section]

\newtheorem{proposition}[theorem]{Proposition}
\newtheorem{lemma}[theorem]{Lemma}

\newtheorem{corollary}[theorem]{Corollary}

{\Alph{theorem}}

\newtheorem*{conjugation}{Standard Forms}
\newtheorem*{class inequality}{Class Inequality}
\newtheorem*{rank inequality}{Rank Inequalities}

\theoremstyle{definition}

\newtheorem{example}[theorem]{Example}

\newtheorem{remark}[theorem]{Remark}
\newtheorem{conjecture}[theorem]{Conjecture}

\newtheorem{chunk}[theorem]{}

\newenvironment{bfchunk}{\begin{chunk}\textit}{\end{chunk}}

\theoremstyle{remark}

\newtheorem{stepp}{Step}

\newtheorem*{Remark}{Remark}

\numberwithin{equation}{theorem}

\newcommand{\btensor}{\boxtimes}

\newcommand{\prclass}[2]{\operatorname{proj\,class}_{#1}#2}
\newcommand{\frclass}[2]{\operatorname{free\,class}_{#1}#2}
\newcommand{\flclass}[2]{\operatorname{flat\,class}_{#1}#2}
\newcommand{\genclass}[2]{\operatorname{{\mathcal P}-class}_{#1}#2}
\newcommand{\cpxx}[2]{{{#1}^{#2}}}
\newcommand{\compx}[2]{{({#1}^{#2})_{\scriptscriptstyle\bullet}}}
\newcommand{\comp}[1]{{{#1}{}_{\scriptscriptstyle\bullet}}}
\newcommand{\filt}[2]{{{#1}^{#2}}}
\newcommand{\Filt}[2]{{\{{#1}^{#2}\}}}

\newcommand{\Comp}[1]{{\mathsf C}(#1)}
\newcommand{\Compplus}[1]{{\mathsf C}_+(#1)}
\newcommand{\Diff}[1]{{\mathsf\Delta}(#1)}
\newcommand{\dfm}[1]{{{#1}{}_{\mathsf{\!\Delta}}}}

\begin{document}

\title[Differential modules]
{Class and rank of differential modules}

\date{\today}

\author[L.~L.~Avramov]{Luchezar L.~Avramov}

\address{Department of Mathematics, University
   of Nebraska, Lincoln, NE 68588, U.S.A.}  \email{avramov@math.unl.edu}

\author[R.-O.~Buchweitz]{Ragnar-Olaf Buchweitz} \address{Department of
Computer and Mathematical Sciences, University of Toronto at Scarborough,
Toronto, ON M1A 1C4, Canada} \email{ragnar@math.toronto.edu}

\author[S.~Iyengar]{Srikanth Iyengar} \address{Department of
Mathematics, University of Nebraska, Lincoln, NE 68588, U.S.A.}
 \email{iyengar@math.unl.edu}

\thanks {Research partly supported by NSF grant DMS 0201904 (L.L.A.),
NSERC grant 3-642-114-80 (R.O.B.), and NSF grant DMS 0442242 (S.I.)}

 \begin{abstract}
A differential module is a module equipped with a square-zero
endomorphism.  This structure underpins complexes of modules over
rings, as well as differential graded modules over graded rings.
We establish lower bounds on the class---a substitute for the length
of a free complex---and on the rank of a differential module in terms
of invariants of its homology.  These results specialize to basic theorems
in commutative algebra and algebraic topology.  One instance
is a common generalization of the equicharacteristic case of the New
Intersection Theorem of Hochster, Peskine, P.\ Roberts, and Szpiro,
concerning complexes over commutative noetherian rings, and of
a theorem of G.\ Carlsson on differential graded modules over
graded polynomial rings.
 \end{abstract}


\subjclass[2000]{Unclassified}

\maketitle

\tableofcontents

\section*{Introduction}

This paper has its roots in a confluence of ideas from commutative algebra
and algebraic topology.  Similarities between two series of results and
conjectures in these fields were discovered and efficiently exploited
by Gunnar Carlsson more than twenty years ago.  On the topological
side they dealt with finite CW complexes admitting free torus actions;
on the algebraic one, with finite free complexes with homology of finite
length. However, no single statement---let alone common proof---covers
even the basic case of modules over polynomial rings.

In this paper we explore the commonality of the earlier results and
prove that broad generalizations hold for all commutative algebras over
fields. They include both Carlsson's theorems on differential graded
modules over graded polynomial rings and the New Intersection Theorem for
local algebras, due to Hochster, Peskine, P.\ Roberts, and Szpiro.  They
also suggest precise statements about matrices over commutative rings,
that imply conjectures on free resolutions, due to Buchsbaum, Eisenbud,
and Horrocks, and conjectures on the structure of complexes with almost
free torus actions, due to Carlsson and Halperin.  These conjectures
are among the fundamental open questions on both sides of this narrative.

The focus here is on a simple construct: a module over an associative
ring $R$, equipped with an $R$-linear endomorphism of square zero. We
call these data a \emph{differential $R$-module}.  They are part of the
structure underlying the familiar and ubiquitous notions of complex or
differential graded module.  Differential modules as such appeared already
five decades ago in Cartan and Eilenberg's treatise \cite{CE}, where
they are assigned mostly didactic functions. Our goal is to establish
that these basic objects are of considerable interest in their own right.

To illustrate the direction and scope of the generality so gained,
take a complex
 \begin{equation*}
P=\quad 0\lra P_l\xra{\ \dd_l\ }P_{l-1}\lra\cdots\lra P_1\xra{\
\dd_1\ }P_{0}\lra0
 \end{equation*}
of finite free modules over a ring $R$.  The module $P=\bigoplus_n
P_n$ with endomorphism $\delta=\bigoplus_n \dd_n$ is a
differential $R$-module $\dfm P$.  With respect to an obvious
choice of basis for the underlying free module, $\delta$ is
represented by a block triangular matrix
 \[
 A=\begin{bmatrix}
0 & A_{01} & 0 & \dots & 0 & 0\\
0 & 0 & A_{12} & \dots & 0 & 0\\
0 & 0 & 0 & \dots & 0 & 0\\
\vdots & \vdots & \vdots & \ddots & \vdots & \vdots\\
0 & 0 & 0 & \dots & 0 & A_{l-1,l}\\
0 & 0 & 0 & \dots & 0 & 0
 \end{bmatrix}
\quad\text{with}\quad A^2=0\,,
\]
Results on finite free complexes are equivalent to statements
about such matrices.

The key contention, supported by our results, is that such
statements should extend in suitable form to \emph{any} strictly
upper triangular matrix
  \[
 A=\begin{bmatrix}
0 & A_{01} & A_{02} & \dots & A_{0\,l-1} & A_{0\,l}\\
0 & 0 & A_{12} & \dots & A_{1\,l-1} & A_{1\,l}\\
0 & 0 & 0 & \dots & A_{2\,l-1} & A_{2\,l}\\
\vdots & \vdots & \vdots & \ddots & \vdots & \vdots\\
0 & 0 & 0 & \dots & 0 & A_{l-1\,l}\\
0 & 0 & 0 & \dots & 0 & 0
 \end{bmatrix}
\quad\text{with}\quad A^2=0\,,
 \]
Matrices of this type arise from sequences of submodules
 \[
\Filt Fn=\quad 0\subseteq\filt F0\subseteq\filt
F1\subseteq\cdots\subseteq\filt Fl=D
 \]
in a differential $R$-module $D$ with endomorphism $\delta$,
satisfying for every $n$ the conditions: $\filt Fn/\filt F{n-1}$
is free of finite rank and $\delta(\filt F{n})\subseteq \filt
F{n-1}$.  We say that $\Filt Fn$ is a \emph{free differential
flag} with $(l+1)$ folds in $D$. When $D$ admits such a flag we
say that its \emph{free class} is at most $l$, and write $\frclass
RD\le l$; else, we set $\frclass RD=\infty$. The \emph{projective
class} of $D$ is defined analogously, and is denoted $\prclass
RD$. Note that if $D$ has finite free (respectively, projective)
class, then it is necessarily finitely generated and free
(respectively, projective).

The \emph{homology} of $D$ is the $R$-module
$\hh{D}=\Ker(\delta)/\Image(\delta)$.  A central result of this
paper links the size of its annihilator, $\Ann_R\hh D$, to the
class of $D$ by a

\begin{class inequality}
Let $R$ be a noetherian commutative algebra over a field and $D$ a
finitely generated differential $R$-module.  One then has
 \[
\prclass RD \geq \height I\quad\text{where}\quad I=\Ann_R\hh D \,.
 \]
\end{class inequality}
The example $D=\dfm K$, where $K$ is the Koszul complex on $d$
elements generating an ideal of height $d$, shows that the
inequality cannot be strengthened in general.
 
For the differential module $\dfm P$ defined above one has 
\[
l\ge \prclass R{\dfm P}\quad\text{and}\quad
\hh{\dfm P}= \bigoplus_n\HH nP\,,
\]
so the New Intersection Theorem follows from the Class Inequality. The hypothesis that $R$
contain a field is due to the use in our proof of Hochster's big Cohen-Macaulay modules
\cite{Ho:CBMS}. The conclusion holds whenever such modules exist, in particular, when
$\dim R\le3$, see \cite{Ho:JA}, or when $R$ is Cohen-Macaulay.  P.~Roberts proved that the
Intersection Theorem holds for all noetherian commutative rings $R$, see \cite{Rb:Camb},
and we conjecture that so does the Class Inequality.

The Class Inequality, a result about commutative rings in general, has its origin in the
study of free actions of the group $(\BZ/2\BZ)^d$ on a CW complex $X$.  Carlsson
\cite{Ca:inv} proved that over a polynomial ring $S=\BF_2[x_1,\dots,x_d]$ a
\emph{differential graded} module $C$ with $\rank{\BF_2}\hh C$ finite and non-zero has
$\frclass SC\geq d$ and used this result to produce obstructions for such actions.  In
\cite{Ca:ln} he conjectured that $\rank SC\ge2^d$ always holds, and showed that a positive
answer implies $\sum_n\rank{\BF_2}{\HH n{X,\BF_2}}\ge2^d$.  This is a counterpart to
Halperin's question as to whether an almost free action of a $d$-dimensional real torus
forces $\sum_n\rank{\BQ}{\HH n{X,\BQ}}\ge2^d$, see \cite{Hl}.

Carlsson \cite{Ca:ams} verified his conjecture for $d\le3$, and Allday
and Puppe \cite{AP} proved that $\rank SC\ge2d$ always holds.  We subsume
these results into the following

\begin{rank inequality}
Let $R$ be a commutative noetherian ring, $D$ a differential
$R$-module of finite free class, and set $d=\height\Ann_R\hh D$.
One then has
 \[
\rank RD\ge\begin{cases}
 2d &\text{when $d\leq 3$ or $R$ is an algebra over a field;}\\
 8 &\text{when $d\geq 3$ and $R$ is a unique factorization domain\,.}
 \end{cases}
 \]
 \end{rank inequality}

We conjecture that an inequality $\rank RD\geq 2^{d}$ always holds.  If it 
does, it will settle the conjectures of Carlsson and Halperin, and will go a 
long way towards confirming a classical conjecture of Buchsbaum, Eisenbud,
and Horrocks: Over a local ring $R$, a \emph{free resolution} $P$ of a
non-zero $R$-module of finite length and finite projective dimension
should satisfy $\rank R{P_n}\ge\binom{\dim R}n$ for each $n$; see
\cite{BE}, \cite{Ha}.

The Rank Inequality is proved in Section \ref{Rank inequalities}.
The Class Inequality is established in Section \ref{Class
inequality. II}, via a version for local rings treated in Section
\ref{Class inequality. I}. The environment for these arguments is
\emph{homological algebra of differential modules}.  References for
the basic formal properties needed to lay out our arguments and a
language fit to express our results are lacking.  In Sections
\ref{Differential modules} and \ref{Differential flags} we close this
gap, guided by the well understood models of complexes and of
differential graded modules.

Such a transfer of technology encounters subtle obstacles.  Complexes and
differential graded modules are endowed with gradings for which the
differential is homogeneous. The use of these gradings in homological
arguments is so instinctive and pervasive that intuition may falter when
they are not available.  

A filtration with adequate properties can sometimes compensate
for the absence of a grading.  This observation led us to the concept
of differential flags.  Their study earned unexpected dividends. One
is an elementary description of matrices that admit a certain standard
form over a local ring. To express it succinctly, let $0_r$ and $1_r$
denote the $r\times r$ zero and identity matrices, respectively.

\begin{conjugation}
\label{standard form} 
Let $A=(a_{ij})$ be a $2r\times 2r$ strictly upper
triangular matrix  with entries in a commutative local
ring $R$.  There is a matrix $U\in{\operatorname {GL}}(2r,R)$ such
that
 \[
U\!AU^{-1}=\begin{bmatrix}0_r & 1_r \\ 0_r & 0_r\end{bmatrix}
 \]
if, and only if, the solutions in $R$ of the linear system of
equations
 \[
\sum_{j=1}^{2r}a_{ij}x_j=0\quad\text{for}\quad i=1,\dots,2r
 \]
are precisely the $R$-linear combinations of the columns of $A$.
 \end{conjugation}

This result appears in Section \ref{Square-zero matrices}.  It is
a consequence of a theorem proved in Section \ref{Differential flags}
for arbitrary associative rings: If $D$ is a differential module
that admits a projective flag, then $\hh D=0$ implies that $D$ is
\emph{contractible}; that is, $D\cong C\oplus C$ with differentiation
the map $(c',c'')\mapsto(c'',0)$.  Such a statement is needed to
get homological algebra started. The corresponding result for bounded
complexes of projective modules holds for trivial reasons.  On the other
hand, not every differential module $D$ with finitely generated projective
underlying module and $\hh D=0$ is contractible. Simple 
examples  are given in Section \ref{Differential modules}.

Our migration to the category of differential modules from the
more familiar environment of complexes was motivated in part by the
investigation in \cite{ABIM} of ``levels'' in derived categories.
The treatment of differential modules presented here is intentionally
lean, concentrating just on what is actually needed to present and prove
the results. A detailed analysis of the homological or homotopical
machinery that differential modules are susceptible of will follow
in \cite{ABCIP}.

\section{Differential modules}
\label{Differential modules}

In this paper $R$ is an associative ring, and rings act on their modules
from the left.  Right $R$-modules are identified with modules over $R\op$,
the opposite ring of $R$.

In this section we provide background on differential modules.  Under the
name `modules with differentiation' the concept goes back to the monograph
of Cartan and Eilenberg \cite{CE}, where it appears twice: in Ch.\ IV, \S\S
1,2 preceding the introduction of complexes, and in Ch.\ XV,
\S\S 1--3 in the construction of spectral sequences.

 \begin{bfchunk}{Differential modules.}
 \label{diff:definitions}
A \emph{differential $R$-module} is an $R$-module $D$ equipped
with an $R$-linear map $\delta^D\col D\to D$, called the
\emph{differentiation} of $D$, satisfying $(\delta^D)^2=0$.
Sometimes we say a pair $(D,\delta)$ is a differential module,
implying $\delta=\delta^D$.

A \emph{morphism} of differential $R$-modules is a homomorphism $\phi\col D\to E$ of
$R$-modules that commutes with the differentiations: $\delta^E\circ\phi =
\phi\circ\delta^D$. If $D$ is a direct summand, as a differential module, of a
differential $R$-module $E$, we say $D$ is a \emph{retract} of $E$; this
distinguishes it from direct summands of the $R$-module $E$.

Note that the category of differential $R$-modules can be identified with the category of
modules over $R[\eps]$, the ring of dual numbers over $R$, see \cite[ Ch.\ IV]{CE}.  In
particular, the following assertions hold: $\Ker(\phi)$ and $\Coker(\phi)$ are differential
$R$-modules; differential $R$-modules and their morphisms form an abelian category, 
denoted $\Diff R$; this category has arbitrary limits and colimits; the formation of
products, coproducts, and filtered colimits are exact operations.
 \end{bfchunk}

\begin{bfchunk}{Homology.}
 \label{difftomod}
For every differential $R$-module $D$ set
 \[
\BB(D) = \Image(\delta^D) \qquad \text{and}\qquad \ZZ(D) = \Ker(\delta^D)\,.
 \]
These submodules of $D$ satisfy $\BB(D)\subseteq \ZZ(D)$ because
$\delta^2=0$. The quotient module
\[
\hh D = \ZZ(D)/\BB(D)
\]
is the \emph{homology} of $D$. We say $D$ is \emph{acyclic} when
$\hh D=0$.

Homology is a functor from $\Diff R$ to the category of $R$-modules.  It commutes with
products, coproducts, and filtered colimits. A \emph{quasi-isomorphism} is a morphism of
differential modules that induces an isomorphism in homology; the symbol $\simeq$
indicates quasi-isomorphisms, while $\cong$ is reserved for isomorphisms.

Every exact sequence of morphisms of differential $R$-modules
 \[
0\to D\to D'\to D''\to0
 \]
induces a \emph{homology exact triangle} of homomorphisms of
$R$-modules
 \begin{equation}
 \label{diff:triangle}
 \begin{gathered}
\xymatrixcolsep{.8pc} \xymatrixrowsep{1.2pc} \xymatrix{
\hh{D}
\ar@{<-}[dr]\ar@{->}[rr] & & \hh {D'} \ar@{->}[dl] \\ & \hh {D''}}
 \end{gathered}
 \end{equation}
For a proof, see \cite[Ch.\ VI, (1.1)]{CE} and the remark following it.

The \emph{suspension} of a differential $R$-module $(D,\delta)$ is
the differential $R$-module
 \begin{equation}
 \label{suspension}
\shift D = (D,-\delta)\,.
 \end{equation}
Suspension is an automorphism of $\Diff R$ of order two; one has
\[
\hh{\shift D}\cong\hh D\,.
\]
Let $\phi\col D\to E$ be a morphism of differential $R$-modules.
The pair
 \begin{equation}
 \label{cone}
\cone\phi=\big(D \oplus E\,,\ (d,e)\mapsto(-\delta^D(d),
\delta^E(e) + \phi(d))\big)
 \end{equation}
is a differential module, called the \emph{cone} of $\phi$.

Obvious morphisms define an exact sequence of differential $R$-modules
 \begin{equation}
   \label{cone:ses} 0\to E \to \cone\phi \to \shift D\to 0\,.
 \end{equation}
Given the isomorphism $\hh D \cong \hh{\shift D}$, it is readily
verified that the homology exact triangle associated with this
exact sequence yields an exact triangle:
 \begin{equation}
\label{cone:triangle}
 \begin{gathered}
 \xymatrixcolsep{1pc} \xymatrixrowsep{1.5pc}
\xymatrix{
\hh{D} \ar@{<-}[dr]\ar@{->}[rr]^{\hh\phi} & & \hh {E} \ar@{->}[dl] \\
& \hh {\cone\phi}}
 \end{gathered}
 \end{equation}
Thus, $\phi$ is a quasi-isomorphism if and only if $\cone\phi$ is acyclic.
 \end{bfchunk}

 \begin{bfchunk}{Compression.}
 \label{complextodiff}
Let $\Comp R$ denote the category of complexes over $R$ with chain maps
of degree $0$ as morphisms.  We display complexes in the form
 \[
X=\quad \cdots\lra X_n\xra{\dd^X_n\ } X_{n-1}\lra\cdots
 \]

The \emph{compression} of a complex $X$ is the differential
$R$-module
 \[
\dfm X=\bigg(\bigoplus_{n\in\BZ} X_n, \bigoplus_{n\in\BZ} \dd^X_n\bigg)\,.
 \]
Compression defines a functor $\Comp R\to\Diff R$. It preserves
exact sequences and quasi-isomorphisms, commutes with colimits,
suspensions, and cones, and satisfies
 \begin{equation*}
 \label{complextodiff:homology}
 \hh{\dfm X} = \bigoplus_{n\in\BZ}\HH nX\,.
 \end{equation*}
\end{bfchunk}

Compression identifies complexes as \emph{graded differential
modules}, but offers no help in studying differential modules.  A
functor in the opposite direction does:

 \begin{bfchunk}{Expansion.}
 \label{difftocomplex}
The \emph{expansion} of a differential $R$-module $D$ is the
complex
 \[
\comp D =\quad\cdots\lra D\xra{\ \delta^D\ }D\xra{\ \delta^D\ }
D\lra\cdots
 \]
Expansion is a functor $\Diff R\to\Comp R$ that commutes with
limits, colimits, suspensions, and cones, and preserves exact
sequences. Since one has
 \begin{equation*}
 \label{difftocomplex:homology}
\HH n{\comp D} = \hh D
 \quad\text{for every}\quad n\in\BZ\,,
\end{equation*}
expansion preserves quasi-isomorphisms as well.
 \end{bfchunk}

 \begin{bfchunk}{Contractibility.}
 \label{contractible}
A differential module is \emph{contractible} if it is isomorphic to
 \begin{equation*}
 \label{contractible:special}
(C\oplus C,\delta)
  \quad\text{with}\quad
\delta(c',c'')=(c'',0)\,.
 \end{equation*}
It is evident that every contractible differential $R$-module is acyclic.
 \end{bfchunk}

In certain cases acyclicity implies contractibility.

 \begin{remark}
 \label{regular}
Assume that $R$ is \emph{regular}, in the sense that every $R$-module
has finite projective dimension.  If $D$ is an acyclic differential
$R$-module, such that the underlying $R$-module $D$ is projective,
then $D$ is contractible.

Indeed, the hypothesis $\hh{\comp D}=0$ yields an exact sequence
of $R$-modules
 \[
0\to \Image(\delta^D) \lra D \xra{\ \delta^D\ } D\lra \cdots\lra D
\xra{\ \delta^D\ } D \lra  \Image(\delta^D) \to 0
 \]
containing $\pd_R\Image(\delta^D)$ copies of $D$.
It follows that $\Image(\delta^D)$ is projective.

Choosing a splitting $\sigma$ of the surjection
$D\to\Image(\delta^D)$, set 
\[
D'=\Ker(\delta^D)\quad\text{and}\quad
D''=\Image(\sigma)\,.
\]
Thus, $D=D'\oplus D''$ and $\delta^D\vert_{D''}$ defines an isomorphism $D''\cong D'$.
 \end{remark}

Examples of non-contractible acyclic differential modules
exist, even with finite free underlying $R$-module, over rings
that are close to being regular.

 \begin{example}
 \label{example:normal}
The ring $R=k[x,y,z]/(x^2 + yz)$, where $k$ is a field and $k[x,y,z]$
a polynomial ring, is a hypersurface and a normal domain with isolated
singularity.  Let $D$ be the differential module with underlying module
$R^2$, defined by the matrix
 \[
A=\begin{bmatrix} x  &  y \\ z  & -x \end{bmatrix}
 \]
with $A^2=0$, see Section \ref{Square-zero matrices}.  Either by direct 
computation or by using Eisenbud's \cite{Ei} technique of matrix 
factorizations, it is easy to check that $D$ is acyclic.
However, $D$ is not contractible: $\Image(\delta^D)\subseteq
(x,y,z)D$ implies $\Image(\delta^D)$ is not a direct summand.
 \end{example}

The next result assesses the gap between acyclicity and
contractibility.

Recall that when $X$ and $Y$ are complexes of $R$-modules $\Hom R{X}{Y}$ denotes the
complex of $\BZ$-modules with
\begin{gather*}
\Hom R{X}{Y}_n = \prod_{i\in\BZ}\Hom R{X_i}{Y_{n+i}}\\
\dd(\vartheta) = \delta^{Y}\vartheta
- (-1)^{|\vartheta|}\vartheta\delta^{X}
\end{gather*}
In particular, in $\Hom R{X}{Y}$ the cycles of degree $0$ are the morphisms of complexes
$X\to Y$, and two cycles are in the same homology class if and only if they are homotopic
chain maps.

  \begin{proposition}
 \label{contractible:homotopy}
For a differential $R$-module $D$ the following are equivalent.
 \begin{enumerate}[\quad\rm(i)]
 \item
$D$ is contractible.
 \item
$\Image(\delta)=\Ker(\delta)$ and the following exact sequence of
$R$-modules is split:
 \[
0\lra \Ker(\delta) \lra D\xra{\ \pi\ }  \Image(\delta)\lra 0\,.
 \]
 \item
$\hh{\Hom R{\comp D}{\comp D}}=0$.
 \end{enumerate}
 \end{proposition}

 \begin{proof}
(i) $\implies$ (iii).  We may assume $D$ has the form
\eqref{contractible:special}.  The maps $\chi_n\col D\to D$ given
by $\chi_n(c',c'')=(0,c')$ for each $n\in\BZ$ then satisfy
$\idmap^{D}=\delta\chi_n+\chi_{n-1}\delta$.  If $\alpha\col{\comp
D}\to{\comp D}$ is a cycle of degree $i$, then
$\delta\alpha_n=(-1)^i\alpha_{n-1}\delta$ holds for all $n$. Thus,
$\chi'_n=\chi_{n+i}\alpha_n\col D\to D$ define a homomorphism
$\chi'\col{\comp D}\to{\comp D}$ such that $\alpha=\dd(\chi').$

(iii) $\implies$ (ii). As $\HH 0{\Hom R{\comp D}{\comp D}}$
vanishes, the identity map $\idmap^{\comp D}$ is homotopic to $0$,
so there are homomorphisms $\chi_n\col D\to D$ of $R$-modules,
satisfying $\idmap^D = \delta\chi_n + \chi_{n-1}\delta$ for each
$n\in\BZ$.  Fix some $n$ and set $\chi=\chi_n$.  One has
 \[
\Ker(\delta)=\delta\chi(\Ker(\delta))+\chi_{n-1}\delta(\Ker(\delta))
=\delta\chi(\Ker(\delta))\subseteq\Image(\delta)\,.
 \]
This implies $\Image(\delta)=\Ker(\delta)$.  The map
$\varepsilon=\delta\chi$ satisfies $\delta=\varepsilon\delta$.
Thus, for $E=\Image(\varepsilon)$ one gets
$E\subseteq\Image(\delta)\subseteq E$, hence $E=\Image(\delta)$.
One also has $\varepsilon^2=\varepsilon$, so for every $e\in E$
the map $\sigma=\chi\vert_{E}\col E\to D$ satisfies
$\delta\sigma(e)=e$, hence the sequence in (ii) splits.

(ii) $\implies$ (i).  The argument at the end of Remark \ref{regular}
applies.
 \end{proof}

 No natural differentiation on tensor products of differential modules commutes with
 expansion, defined in \ref{difftocomplex}.  The absence of tensor products on the category of
 differential modules seriously limits the applicability of standard technology.  A more
 frugal structure, defined below, provides a partial remedy to that situation.

 \begin{bfchunk}{Tensor products.}
 \label{diff:tensors}
Let $R'$ be an associative ring and $X$ a complex of $R'$-$R\op$
bimodules.  The \emph{tensor product} of $X$ and $D\in\Diff R$ is
the differential $R'$-module
 \begin{equation*}
  \label{tensor}
X\btensor_R D = \bigg(\bigoplus_{n\in\BZ}(X_n\otimes_RD)\,,\ x\otimes
d \mapsto \dd^{X}(x) \otimes d + (-1)^{|x|}x\otimes
 \delta^D(d)\bigg)
 \end{equation*}
where $|x|$ denotes the degree of $x$.  Tensor product defines a
functor
 \[
-\btensor_R-\col\Comp{R'\otimes_{\BZ}R\op}\times\Diff R\lra\Diff{R'}\,.
 \]

Whenever needed, modules or bimodules are considered to be
complexes concentrated in degree $0$.  Thus, $M\btensor_RD$ is
defined for every $R'$-$R\op$-bimodule $M$; as it is equal to
$(M\otimes_RD,M\otimes_R\delta)$, we sometimes write $M\otimes_RD$
in place of $M\btensor_RD$.

Tensor products commute with colimits on both sides. We
collect some further properties, using equalities to denote
canonical isomorphisms. For every complex $X$ of
$R'$-$R\op$-bimodules and every differential $R$-module $D$, one
has
 \begin{alignat}{2}
 \label{btensor:expansion}
\comp{(X\btensor_RD)} &= X\otimes_R \comp D
&&\qquad\text{in}\quad\Comp{R'}
 \\
 \label{btensor:homology}
\hh{X\btensor_RD}&=\HH n{X\otimes_R \comp D} &&\qquad\text{for
each}\quad n\in\BZ
\\
 \intertext{For every complex $W$ of $R''$-$R'\op$ bimodules one has}
 \label{btensor:associative}
(W\otimes_{R'}X)\btensor_RD &=W\btensor_{R'}(X\btensor_RD)
&&\qquad\text{in}\quad \Diff{R''}
\\
 \label{btensor:compression}
W\btensor_{R'}\dfm{X}&= \dfm{(W\otimes_{R'}X)}
&&\qquad\text{in}\quad \Diff{R''}
 \\
 \intertext{For the differential $R$-module $R=(R,0)$ one has}
 \label{btensor:rightunit}
X\btensor_RR&=\dfm{X} &&\qquad\text{in}\quad \Diff{R'}
 \\
 \label{btensor:leftunit}
R\btensor_RD&= D
&&\qquad\text{in}\quad \Diff R
\\
 \intertext{For all morphisms $\vartheta\in\Comp{R'\otimes_{\BZ}R\op}$ and
$\phi\in\Diff R$ one has}
 \label{btensor:compcones}
\cone{\vartheta\btensor_RD}&=\cone{\vartheta}\btensor_RD
&&\qquad\text{in}\quad\Diff{R'}
 \\
  \label{btensor:dfmcones}
\cone{X\btensor_R\phi}&\cong X\btensor_R\cone\phi
&&\qquad\text{in}\quad\Diff{R'}
 \end{alignat}
 \end{bfchunk}

 We need to track exactness properties of tensor products.  Evidently, if $D$ is
 contractible, for each complex $X$ of $R'$-$R\op$-bimodules the differential $R'$-module
 $X\btensor_RD$ also is contractible, and hence acyclic.  The next result shows, in
 particular, that $(X\btensor_R-)$ preserves acyclicity under the expected hypotheses on
 $X$.

 \begin{proposition}
 \label{diff:quisms}
Let $X$ and $Y$ be bounded below complexes of
$R'$-$R\op$-bimodules, such that the $R\op$-modules $X_i$ and
$Y_i$ are flat for all $i\in\BZ$.
 \begin{enumerate}[\quad\rm(1)]
 \item The following functor preserves exact sequences and quasi-isomorphisms
 \[
(X\btensor_R-)\col \Diff R\lra \Diff{R'}\,.
 \]
 \item
A quasi-isomorphism $\vartheta\col X\to Y$ in $\Comp{R'\otimes_\BZ
R\op}$ induces for each differential $R$-module $D$ a
quasi-isomorphism of differential $R'$-modules
 \[
\vartheta\btensor_RD\col X\btensor_RD\lra Y\btensor_RD\,.
 \]
 \end{enumerate}
 \end{proposition}

 \begin{proof}
(1) The functor $(X\btensor_R-)$ preserves exact sequences because
the complex $X$ consists of flat $R\op$-modules.  On the other
hand, a morphism $\phi$ of differential modules is a
quasi-isomorphism if and only if its cone is acyclic, see
\eqref{cone:triangle}.  In view of the isomorphism
\eqref{btensor:dfmcones}, to finish the proof it suffices to show
that if a differential module $D$ is acyclic, then so is
${X\btensor_RD}$.

For each integer $n$ define a subcomplex of $X$ as follows:
 \[
X_{\les n}=\quad \cdots\lra0\lra0\lra X_n\lra X_{n-1}\lra\cdots
 \]
It fits into an exact sequence of complexes of right $R$-modules
 \[
0\lra X_{\les n-1}\lra X_{\les n}\lra \shift^n X_n\lra 0\,.
 \]
 This sequence induces an exact sequence of differential modules
 \[
\xymatrixcolsep{1.5pc} \xymatrix{
 0 \ar@{->}[r]
 & X_{\les n-1}\btensor_RD\ar@{->}[r]
 & X_{\les n}\btensor_RD\ar@{->}[r]
 & (\shift^n X_n)\btensor_RD\ar@{->}[r]
 & 0\,.}
 \]
As the $R\op$-module $X_n$ is flat, one has
\[
\hh{(\shift^nX_n)\btensor_RD}\cong X_n\otimes_R\hh D=0\,.
\]
Since $X_{\les n}=0$ holds for $n\ll0$, by induction we may assume $\hh{X_{\les
    n-1}\btensor_RD}=0$ holds as well.  The exact sequence above yields $\hh{X_{\les
    n}\btensor_RD}=0$, so using the equality $X=\bigcup_{n}X_{\les n}$ and the exactness
of colimits one obtains
 \[
 \hh{X\btensor_RD}
 =\hh{\colim_n(X_{\les n}\btensor_RD)}
 =\colim_n\hh{X_{\les n}\btensor_RD}=0\,.
 \]

(2) Note that $W=\cone{\vartheta}$ is a bounded below complex of
flat $R\op$-modules with $\hh{W}=0$.  Arguing as in (1), one sees
it suffices to prove $W\btensor_RD$ is acyclic. By
\eqref{btensor:homology} this is equivalent to proving that the
complex $W\otimes_R\comp D$ is acyclic.

Set $(-)^\vee=\Hom{\BZ}{-}{\BQ/\BZ}$.  As the $\BZ$-module $\BQ/\BZ$ is faithfully
injective, it suffices to prove ${\hh{ W\otimes_R\comp D}}^\vee=0$.  The exactness of
$(-)^\vee$ yields $\hh{{W}^{\vee}}\cong{\hh{W}}^\vee=0$.  It also implies that each
$(W_i)^\vee$ is an injective $R$-module, because $W_i$ is a flat $R\op$-module. Thus,
${W}^{\vee}$ is acyclic, bounded above complex of injective $R$-modules, and so it is
contractible.  This explains the equality in the sequence
 \[
{\hh{W\otimes_R\comp D}}^{\vee}\cong
 \hh{({W\otimes_R\comp D})^{\vee}}\cong
  \hh{\Hom{R}{\comp D}{W^{\vee}}}=0\,.
   \]
Exactness of $(-)^{\vee}$ yields the first isomorphism and
adjointness the second.
 \end{proof}

Finding properties of differential modules that guarantee
exactness of tensor products is a more delicate matter.  It is
discussed in the next section.

\section{Differential flags}
\label{Differential flags}

Throughout this section $R$ denotes an associative ring.  We
introduce and study classes of differential $R$-modules that
conform to classical homological intuition.

\begin{bfchunk}{Flags.}
\label{flags:definitions} Let $F$ be a differential $R$-module.  A
\emph{differential flag} in $F$ is a family $\Filt Fn$ of
$R$-submodules satisfying the conditions
 \begin{align}
 \label{flag:union}
\filt F{n}&\subseteq\filt F{n+1}\,,
 \quad \filt F{-1}=0\,,
 \quad \bigcup_{i\in\BZ}\filt Fi=F\,,
\quad\text{and}\quad
 \\
  \label{flag:nilpotence}
\delta(\filt Fn)&\subseteq \filt F{n-1}
 \end{align}
 for each $n\in\BZ$. Condition \eqref{flag:nilpotence} implies $\filt Fn$ is a
 differential submodule of $F$ and 
 \begin{equation}
 \label{flag:fold}
F_n=\filt Fn/\filt F{n-1}
 \end{equation}
is a differential module with trivial differentiation; we call it
the $n$th \emph{fold} of $\Filt Fn$.  The flag defines for each
$n\in\BZ$ an exact sequence of differential $R$-modules
 \begin{equation}
 \label{flag:sequence}
0\lra \filt F{n-1}\xra{\ \filt{\iota}{n-1}\ }\filt Fn\lra
F_n\lra0\,.
\end{equation}
 \end{bfchunk}

Properties of the folds of a flag affect the character
of a differential module to a stronger degree than do
properties of the underlying module.

\begin{bfchunk}{Types of flags.}
\label{flags:types} A flag $\Filt Fn$ in $F$ is \emph{free}
(respectively, \emph{projective}, \emph{flat}) if every fold $F_n$
has the corresponding property; these conditions are progressively
weaker.

When $\Filt Fn$ is a projective flag one has
$F\cong\bigoplus_{n=0}^\infty F_n$ as $R$-modules.

When $\Filt Fn$ is a flat flag \eqref{flag:sequence} induces an
exact sequence of differential modules
 \begin{equation}
 \label{flag:flat}
\xymatrixrowsep{1.5pc} \xymatrixcolsep{2pc} \xymatrix{
0\ar@{->}[r] & X\btensor_R \filt F{n-1} \ar@{->}[r] & X\btensor_R
\filt F{n}
  \ar@{->}[r] &X\btensor_R F_n\ar@{->}[r] & 0}
 \end{equation}
 for every $X\in\Comp{R\op}$.  Its homology exact triangle \eqref{diff:triangle} has the
 form
 \begin{equation}
 \label{flag:triangle}
 \begin{gathered}
\xymatrixcolsep{1.3pc}
\xymatrixrowsep{2.5pc}
\xymatrix{
 \hh{{X}\btensor_R\filt F{n-1}}
 \ar@{->}[rr]^{\hh{{X}\btensor_R\iota^{n-1}}}
&&\hh{X\btensor_R\filt F{n}}
 \ar@{->}[dl]\\
&\hh{X}\btensor_R F_n
 \ar@{->}[ul]
 }
 \end{gathered}
 \end{equation}
because the differential module $F_n$ is flat and has zero differentiation.
\end{bfchunk}

 \begin{theorem}
 \label{flags:triviality}
Let $D$ be a retract of a differential $R$-module $F$ that admits a
projective flag $\Filt Fn$.

If $D$ is acyclic, then $D$ is contractible.
  \end{theorem}

 \begin{proof}
By Proposition \ref{contractible:homotopy}, it suffices to show
that $\Hom R{\comp D}{\comp D}$ is acyclic. It is a 
retract of $\Hom R{\comp F}{\comp D}$, so we prove that $\hh X=0$
implies $\hh{\Hom R{\comp F}{X}}=0$. 

First we show that $\Hom R{\compx Fn}{X}$ is acyclic by induction
on $n$. Each sequence \eqref{flag:sequence} yields an exact
sequence of complexes of $R$-modules
 \[
0\lra \compx F{n-1}\xra{\ \compx{\iota}{n-1}\ }\compx Fn\lra
\comp{(F_n)}\lra0\,.
 \]
Since $\cpxx Fn=0$ for $n<0$, we may assume $\hh{\Hom R{\compx
F{n-1}}{X}}=0$ for some $n\ge0$.  Each $\iota^{n-1}_i$ is split,
so $\pi^n=\Hom R{\compx\iota{n-1}}{X}$ is surjective, hence the
sequence
 \[
0\to \Hom R{\comp{(F_n)}}{X} \to
     \Hom R{\compx F{n}}{X}\xra{\ \pi^n}
     \Hom R{\compx F{n-1}}{X}\to 0
 \]
is exact.  The complex $\comp{(F_n)}$ has zero differential and
$F_n$ is projective, so one obtains
 \begin{align*}
\hh{\Hom R{\comp{(F_n)}}{X}}&=\hH\bigg(\prod_{i\in\BZ}\shift^i\Hom R{F_n}{X}\bigg) \\
&=\prod_{i\in\BZ}\shift^i\hh{\Hom R{F_n}{X}} \\
&=0\,.
 \end{align*}
The exact sequence and the induction hypothesis yield 
\[
\hh{\Hom R{\compx Fn}{ X}}=0\,.
\]
The first isomorphism below comes from \eqref{flag:union}, the
second is standard:
 \begin{align*}
\hh{\Hom R{\comp F}{X}} &= \hh{\Hom R{\colim_n\compx F{n}}{X}}\\
                        &=\hh{\lim{\!}_n\Hom R{\compx F{n}}{X}}\,.
\end{align*}
Since the limit is taken over the surjective morphisms $\pi^n$ and each
complex in the inverse system is acyclic, the limit complex is  acyclic.
 \end{proof}

The following result complements Proposition \ref{diff:quisms}.

Recall that $\Comp{R'}\otimes_\BZ R\op)$ is the category of
complexes of $R'$-$R\op$-bimodules; let $\Compplus{R'\otimes_\BZ
R\op}$ be its full subcategory consisting of bounded below
complexes.

 \begin{proposition}
 \label{flags:semi-flat}
Let $D$ and $E$ be retracts of differential $R$-modules with flat
flags.
 \begin{enumerate}[\quad\rm(1)]
 \item
The following functor preserves exact sequences and
quasi-isomorphisms:
 \[
(-\btensor_RD)\col \Comp{R'}\otimes_\BZ R\op) \lra \Diff{R'}\,.
 \]
 \item
A quasi-isomorphism $\phi\col D\to E$ of differential $R$-modules
induces for each $X\in\Compplus{R'\otimes_\BZ R\op}$
a quasi-isomorphism of differential $R'$-modules
 \[
X\btensor_R\phi\col X\btensor_RD\lra X\btensor_RE\,.
 \]
 \end{enumerate}
 \end{proposition}

 \begin{proof}
(1) The functor $-\btensor_RD$ preserves exact sequences because
the $R$-module $D$ is flat.  It remains to verify that
$-\btensor_RD$ preserves quasi-isomorphisms.  As $D$ is a retract
of a differential module $F$ with a flat flag, it suffices to
prove that if $\vartheta$ is a quasi-isomorphism of complexes,
then $\vartheta\btensor_RF$ is one of differential modules.

Recall that a morphism of complexes or of differential modules is
a quasi-isomorphism if and only if its cone is acyclic. Therefore,
in view of the isomorphism $\cone{\vartheta\btensor_RF}
=\cone{\vartheta}\btensor_RF$, see \eqref{btensor:compcones}, it
suffices to prove that for each acyclic complex $X$, the
differential module $X\btensor_RF$ is acyclic.

Let $\Filt Fn$ be a flat flag in $F$.  We prove $\hh{X\btensor_R
\filt F{n}}=0$ by induction on $n$. This is obvious for $n<0$, as
one then has $\filt F{n}=0$. If $X\btensor_R\filt F{n-1}$ is
acyclic for some $n\ge0$, then the exact triangle
\eqref{flag:triangle} yields $X\btensor_R \filt F{n}$ is acyclic,
as desired.  The flag $\Filt F{n}$ in $F$ induces a flag
$\Filt{X\btensor_RF}{n}$ in $X\btensor_RF$, so one gets
 \begin{align*}
 \hh{X\btensor_RF}
&=\hh{\colim_n(X\btensor_R\filt F{n})} \\
&=\colim_n\hh{X\btensor_R\filt F{n}}\\
&=\colim_n 0\\
&=0\,.
\end{align*}

(2) Choose a quasi-isomorphism $\rho\col W \to X$, where $W$
is a bounded below complex of flat $R\op$-modules.  In the
commutative diagram
 \[
 \xymatrixrowsep{2.5pc}
 \xymatrixcolsep{3.5pc}
\xymatrix{
 W\btensor_RF
 \ar@{->}[r]^-{W\btensor_R\phi}
 \ar@{->}[d]_-{\rho\btensor_RD}
 &W\btensor_RE
 \ar@{->}[d]^-{\rho\btensor_RE}
\\
 X\btensor_RD
 \ar@{->}[r]^-{X\btensor_R\phi}
 &X\btensor_RE
}
 \]
both vertical maps are quasi-isomorphisms by (1), and
$W\btensor_R\phi$ is a quasi-iso\-mor\-phism by Proposition
\ref{diff:quisms}(1).  Thus, $X\btensor_R\phi$ is a
quasi-isomorphism, as desired.
 \end{proof}

We show by example that the the flag structure of $F$ is
essential for the validity of the preceding theorem. In fact, its
conclusion may fail even if the differential module involved is
free as an $R$-module.

 \begin{example}
 \label{dold2}
Set $R=\BZ/(4)$.  A projective resolution of $k=R/(2)$ is given by
 \[
X = \quad \cdots\lra R\xra{\ 2\ }R\xra{\ 2\ }R\lra0\lra0\lra\cdots
 \]
 The differential $R$-module $D=(R,2\cdot\idmap^R)$ has $\hh D=0$, therefore
 $\hh{X\btensor_RD}=0$, see Proposition \ref{diff:quisms}.(1).  However,
 $\hh{k\btensor_RD}=k$, so the map $\vartheta\btensor_RD$ induced by the augmentation
 $\vartheta\col X\to k$ is not a quasi-isomorphism.
\end{example}

Each flag in a differential module naturally gives rise to a spectral sequence, see
\cite[Ch.\ XV, \S\S 1--3]{CE}. It is used to prove Theorem \ref{horrocks:ufd}.

\begin{bfchunk}{Spectral sequences.}
 \label{spectral sequence}
Let $\Filt F{n}$ be a flag in a differential $R$-module $F$. For each $r\ge1$ the $r$th
page of the spectral sequence is a family $\Gr {}r{}Fn$ of $r$ complexes
 \[
\eGr pr{} =  \quad
\cdots \lra \eGr {}r{i+r}\xra{\ \dd_{i,i+r}\ } \eGr {}r{i}
\xra{\ \dd_{i-r,i}\ } \eGr {}r{i-r}\lra\cdots
\]
of $R$-modules, where $p=0,1,\dots,r-1$ and $i\equiv p \mod r$. The first page is
\[
\eGr 01i = F_i \quad \text{and} \quad \dd_{i-1,i}(x + \filt F{i-1}) = \delta(x) + \filt F{i-2}\,.
\]
Successive pages of the spectral sequence are linked by equalities
\[
\eGr{}{r+1}i = \Ker(\dd_{i-r,i})/\Image(\dd_{i,i+r})
\quad\text{for each pair}\quad (r,i)\in\BN\times\BN\,.
\]
Evidently when $r\geq i+1$ one has $\dd_{i-r,i}=0$ and there is a surjective system
\[
\eGr{}{r}i\lra \eGr{}{r+1}i\lra \eGr{}{r+2}i\lra\cdots
\]
One sets $\eGr{}{\infty}i =\colim_{r} \eGr{}{r}i$.  For each integer $i$, let
\[
\filt{\hh{F}}i=\Image(\hh{\filt Fi}\to\hh F)\,,
\]
where the arrow is induced by the inclusion $\filt Fi\subseteq F$.  The spectral sequence
strongly converges to $\hh F$, in the sense that there are isomorphisms of $R$-modules
 \begin{equation}
\label{spectral:convergence}
\begin{split}
\filt{\hh{F}}i/\filt{\hh{F}}{i-1}\cong\eGr{}{\infty}i
 \quad\text{for each $i\geq 0$}\\
\filt{\hh{F}}{-1}=0 \quad\text{and}\quad \bigcup_{i\in\BZ}\filt{\hh{F}}i=\hh F\,.
\end{split}
 \end{equation}

If $\filt Fl = F$ for some $l\ge0$, then one has $\eGr{}{\infty}i=0$ for
$i\notin[0,l]$,  hence
\begin{equation}
\label{spectral:stablepage}
\eGr{}{\infty}i=\eGr{}{r}i \quad \text{for} \quad r\ge \underset{0\les i\les l}{\max}\{i,l-i\}+1\,.
\end{equation}
 \end{bfchunk}

 Convergence of the spectral sequence above transfers information from $\gr{r}Fn$ to $\hh
 F$. For instance, if the $R$-module $\bigoplus_i\HH i{\gr rFn}$ has finite length for
 some $r$, then so does the $R$-module $\hh F$.  However, this fact does not reduce the
 study of differential modules to that of complexes. The reason is that properties of $\hh
 F$, the primary invariant of $F$, rarely translate into usable information about the
 pages of the spectral sequence. An explicit example is given next.

\begin{example}
  Set $R=k[x,y]/(x^2,xy)$.  By Section \ref{Square-zero matrices} the matrix 
 \[
A=\left[\begin{gathered}
 \xymatrixrowsep{.2pc}
 \xymatrixcolsep{.2pc}
\xymatrix{
0 & x & y & 0\\
0 & 0 & x & y \\
0 & 0 & 0 & 0 \\
0 & 0 & 0 & 0 }\end{gathered}\right]
\]
defines a differential module $F$. The complex $\Gr {}1{}Fn$ has the form
 \[
\cdots\lra0\lra R^2\xra{\ [x\ y]\ }R\xra{\ \ x\ \ }R\lra0\lra\cdots
 \]
The length of $\hh F$ is finite, because $\fp=(x)$ is the only
non-maximal ideal of $R$, the local ring $R_\fp$ is equal to the
field $k(y)$ and $\rank {k(y)}(A_\fp)=2$.  On the other hand, one
has $\Gr 020Fn\cong k[y]$, and this module has infinite
length.
 \end{example}

 Next we define invariants that are central to this paper.  The terminology is modelled on the
 usage of `class' in group theory to measure the shortest length of a filtration with
 subquotients of a certain type, such as in a `nilpotent group of class $l$'.  The length
 of a `solvable' free differential graded modules over graded polynomial ring, introduced
 in \cite[Def.\ 9]{Ca:inv}, is related to its free class, defined below.

 \begin{bfchunk}{Class.}
 \label{class}
We define the \emph{flat class} of a differential $R$-module $F$
to be the number
 \[
\flclass RF= \inf\left\{l\in\BN\left|
   \begin{gathered}
 \text{$F$ admits a flat flag}\\
 \text{$\Filt Fn$ with }\filt Fl=F
   \end{gathered}
\right\}\right.
 \]
The \emph{projective class} of $F$ over $R$, denoted  $\prclass
RM$, and its \emph{free class}, denoted $\frclass RF$, are defined
similarly.  We list some simple properties of these invariants. In
statements valid for any flavor of the definition, we let
$\mathcal P$- stand for either `flat', `projective', or `free',
and let $\genclass RM$ denote the corresponding number.
 \begin{enumerate}[\rm(1)]
 \item
$\genclass RF=\infty$ if and only if $F$ admits no finite $\mathcal
P$-flag.
 \item
$\genclass RF=0$ if and only if $F$ is a $\mathcal P$-module and
$\delta^F=0$.
 \item
If $F$ is a contractible, non-zero, projective (respectively,
flat) module over $R$, then $\prclass RF=1$ (respectively,
$\flclass RF=1$).
 \item
If $F'$ and $F''$ are differential $R$-modules, then flat and
projective class satisfy
 \[
\genclass R{(F'\oplus F'')}=\max\{\genclass RF'\,,\,\genclass
RF''\}\,.
 \]
\item
If $0\to F\to F'\to F''\to 0$ is an exact sequence of differential
$R$-modules, then
  \[
\genclass R{F'}\le \genclass R{F}+\genclass R{F''}+1\,.
  \]
 \item
If $F=\dfm P$, where $P$ is a non-zero, bounded below complex of
$\mathcal P$-modules, then
 \[
\genclass RF\le\card\{i\in\BZ\mid P_i\ne0\}-1\,.
 \]
 \item
For every $F$ the following inequalities hold:
 \[
\flclass RF\le\prclass RF\le\frclass RF\,.
 \]
 \item
When $F$ is finitely generated and $R$ is noetherian, one has
 \[
\flclass RF=\prclass RF\,.
 \]
\item
If $R$ is an IBN ring, see Section \ref{Square-zero matrices}, and $F$ has a free flag, then 
 \[
\frclass RF\le\rank RF\,.
 \]
 \item
If $R\to S$ is a homomorphism of rings, then 
 \[
\genclass S{(S\otimes_RF)}\le\genclass RF\,.
 \]
\end{enumerate}

Indeed, (6) follows from (5) and (2).  For (8), note that if $\Filt Fn$
is a flat flag in $F$, then each fold $F_n$ is finitely presented, hence
it is projective.  Ranks of free modules, when defined, are additive in
exact sequences: this gives (9).  The other assertions follow directly
from the definition of class.
 \end{bfchunk}

A moment of reflection shows why non-trivial lower bounds on $\genclass
RF$ may be hard to obtain.  One method for obtaining such bounds is
given by the following technical result, distilled from the proof of
\cite[Thm.\ 16]{Ca:inv}.

 \begin{proposition}
 \label{carlsson}
Consider a sequence of complexes
 \[
\cpxx X{s}\xra {\cpxx{\vartheta}{s}} \cpxx X{s-1} \lra \cdots \lra
\cpxx X{0} \xra {\cpxx{\vartheta}{0}} X
 \]
 of $R\op$-modules with the following property:
 \begin{equation}
 \tag{a}
\hh{\cpxx{\vartheta}{n}}=0 \quad\text{for}\quad n=0,1,\dots,s\,.
 \end{equation}
If $\pi\col F\to D$ is a morphism  of differential $R$-modules
satisfying the condition
 \begin{equation}
 \tag{b}
\hh{{\vartheta}\btensor_R\pi}\ne 0\quad\text{for}\quad
{\vartheta}=\cpxx{\vartheta}{0}\circ\cdots\circ\cpxx{\vartheta}{s}\,,
 \end{equation}
then the  inequality below  holds:
 \[
\flclass RF\ge s+1\,.
 \]
 \end{proposition}

 \begin{proof}
Let $\Filt F{n}$ be a flat flag with $\filt F{s}=F$, and let
$\iota^{n}\col \filt F{n}\to \filt F{n+1}$ denote the inclusions of differential 
submodules. Form the morphisms of complexes 
 \begin{align*}
 \cpxx{\theta}n=\cpxx{\vartheta}{0}\circ\cdots\circ\cpxx{\vartheta}{n}
 &\col X^n \lra X  \quad\text{for }n=0,\dots,s\,;
\\
 \cpxx{\theta}{-1}=\idmap^X &\col X \lra X\,,
 \end{align*}
and the morphisms  of differential modules
 \[
 \phi^n=\pi\circ\iota^{s}\circ\cdots\circ\iota^{n}
 \col\filt F{n} \lra D
 \quad\text{for }n=-1,\dots,s\,.
 \]
By descending induction on $n$ we will show that the map
 \begin{equation*}
\hh{\cpxx{\theta}n\btensor_R\phi^n} \col \hh{\cpxx
X{n}\btensor_R\filt F{n}} \to \hh{X\btensor_RD}
 \end{equation*}
is non-zero for each integer $n\in[-1,s]$.  This contradicts $\filt F{-1}=0$.

One has $F^s=F$, so (b) is the desired assertion for $n=s$. Assume
that it holds for some $n\in[0,s]$. The exact triangle
\eqref{flag:triangle} yields a commutative ladder 
 \[
\xymatrixrowsep{2.6pc}
\xymatrixcolsep{5pc}
\xymatrix{
 \hh{\cpxx X{n}\btensor_R\filt F{n-1}}
 \ar@{->}[r]_-{\hh{\cpxx{\vartheta}{n}\btensor_R \filt F{n-1}}}
 \ar@{->}[d]_{\hh{\cpxx X{n}\btensor_R\iota^{n-1}}}
&\hh{\cpxx X{n-1}\btensor_R\filt F{n-1}}
 \ar@{->}[d]^{\hh{\cpxx X{n-1}\btensor_R\iota^{n-1}}} \\ 
\hh{\cpxx X{n}\btensor_R\filt F{n}}
 \ar@{->}[r]_{\hh{\cpxx{\vartheta}{n}\btensor\filt F{n}}}
 \ar@{->}[d]
&\hh{\cpxx X{n-1}\btensor_R\filt F{n}}  \ar@{->}[d]\\
\hh{\cpxx X{n}}\btensor_R F_n
 \ar@{->}[r]_{\hh{\cpxx{\vartheta}{n}}\btensor_RF_n}
&\hh{\cpxx X{n-1}}\btensor_RF_n}
 \]
with exact rows.  As $\hh{\cpxx{\vartheta}{n}}\btensor_RF_n=0$ holds by condition
(a), one has an inclusion
 \[
 \Image\hh{\cpxx X{n-1}\btensor_R\iota^{n-1}}
\supseteq
 \Image{\hh{\cpxx{\vartheta}{n}\btensor_R\filt F{n}}}\,.
 \]
In view of the definitions of $\theta^n$ and $\phi^n$ and of the
induction hypothesis, it yields
 \begin{align*}
\Image\hh{\cpxx{\theta}{n-1}\btensor_R\phi^{n-1}}
 &=\hh{\cpxx{\theta}{n-1}\btensor_R\phi^{n}}
 (\Image\hh{\cpxx X{n-1}\btensor_R\iota^{n-1} })
\\
 &\supseteq\hh{\cpxx{\theta}{n-1}\btensor_R\phi^{n}}
 (\Image{\hh{\cpxx{\vartheta}{n}\btensor_R\filt F{n}}})
\\
 &=\Image\hh{\cpxx{\theta}{n}\btensor_R\phi^{n}}
\\
 &\ne0\,.
 \end{align*}

The induction step is now complete, and so the proposition is
proved.
 \end{proof}

\section{Class inequality. I}
\label{Class inequality. I}

For most of this section $(R,\fm,k)$ is a \emph{local ring},
meaning that $R$ is commutative and noetherian, $\fm$ is its
unique maximal ideal, and $k=R/\fm$ its residue field.

The next theorem is the main step towards establishing the Class
Inequality announced in the introduction.  It is in some respects sharper
than the global version, see Theorem \ref{intersection}. The proof is
given at the end of the section.

 \begin{theorem}
 \label{intersection:local}
Let $(R,\fm,k)$ be a local ring, $F$ a differential $R$-module,
and $D$ a retract of $F$ such that the $R$-module $\hh D$ has
non-zero finite length.

When $R$ has a big Cohen-Macaulay module one has:
\[
  \flclass RF\geq \dim R\,.
\]
 \end{theorem}

 We pause to recall existence results for big Cohen-Macaulay modules, and to discuss
 antecedents of the theorem.

\begin{bfchunk}{Big Cohen--Macaulay modules.}
 \label{bigCM}
Recall that an $R$-module $M$ is \emph{big Cohen-Macaulay} if for some
system of parameters $\bsx=x_1,\dots,x_d$ of $R$ the element $x_{i}$
is not a zero divisor on $M/(x_1,\dots,x_{i-1})M$ for $i=1,\dots,d$
and $M\ne(\bsx)M$.  It is not known whether every $R$ has such a module,
but many important cases are covered:

Big Cohen-Macaulay modules exist when $R$ contains a field as a subring,
due to a celebrated construction of Hochster \cite{Ho:CBMS}. They exist
also over all local rings of dimension at most $3$: the difficult case of
dimension $3$ is settled by Hochster \cite{Ho:JA} using Heitmann's proof
of the direct summand conjecture in dimension $3$. Any Cohen-Macaulay
ring $R$ is a big Cohen-Macaulay module over itself.
 \end{bfchunk}

We recall a fundamental result in commutative algebra:

 \begin{remark}
 \label{NIT}
The New Intersection Theorem reads:  Let $R$ be a local ring and let
 \[
P=\quad\cdots\lra0\lra P_l\lra\cdots\lra P_0\lra0\lra\cdots
 \]
be a complex of finite free modules with $P_0\ne 0\ne P_l$; if $P$ is not
exact, and $\length(\HH nP)$ is finite for each $n$, then $l\ge\dim R$.

Hochster, Peskine and Szpiro, and P.\ Roberts established the
theorem when $R$ has a big Cohen-Macaulay module.  {\v C}ech
complexes play a role in all these proofs, cf.\ the discussion in
\cite[pp.\ 82--86]{Rb:Sem}.  In mixed characteristic the
theorem was proved by Roberts, using local Chern classes. His
monograph \cite{Rb:Camb} contains detailed arguments and develops
the necessary intersection theory. The technology powering this
portion of the proof has no analog for differential modules at
present.
 \end{remark}

Theorem \ref{intersection:local} contains the New Intersection
Theorem for rings with big Cohen-Macaulay modules, see Remark
\ref{class}(6), and vastly generalizes the next result.

 \begin{remark}
 \label{CT}
Carlsson's theorem, \cite[Thm.\ 16]{Ca:inv}, may be stated as follows:
Let $R$ be a {polynomial ring} in $d$ variables of {positive degree} over
a field $k$ and $F$ a differential module with {finitely generated graded
free} underlying module and $\delta^F$ {homogeneous of degree $-1$};
every {homogeneous} free flag $\Filt Fn$ in $F$ then has $\filt Fd\ne F$.
 \end{remark}

To prove Theorem \ref{intersection:local} we transplant an idea
from \cite{Ca:inv}, see Proposition \ref{carlsson}, utilize {\v C}ech
complexes, see  \ref{stables}, and introduce two novel ingredients.

One is the determination of a framework for stating and proving a common
generalization of the two theorems above; it is given  by differential
modules with flat flags.  Their properties are put to full use:
almost every result established in Sections \ref{Differential modules}
and \ref{Differential flags} participates in the proofs of Theorems
\ref{intersection:local} and \ref{trivialflags:local}.

A second new ingredient is Theorem \ref{trivialflags:local} below, a homological version
of Nakaya\-ma's Lemma. A similar statement for bounded below complexes follows easily by
inspecting the augmentation map to the non-vanishing homology module of lowest degree. In
the ungraded world of differential modules such a map simply does not exist, so a
completely new approach is needed.  To this end we adapt a \emph{d\'evissage} procedure
introduced by Dwyer, Greenlees, and Iyengar \cite[\S5]{DGI}.

 \begin{theorem}
 \label{trivialflags:local}
Let $(R,\fm,k)$ be a local ring and $M$ an $R$-module with $\fm M\ne M$.
Let $D$ be a retract of a differential $R$-module $F$ that admits a
flat flag.

If $\hh D$ is finitely generated and $M\otimes_RD$ is acyclic, then
$D$ is acyclic.

If, in addition, $F$ admits a projective flag, then $D$ is contractible.
 \end{theorem}

\begin{proof}
The second assertion follows from the first one and Theorem
\ref{flags:triviality}.

We prove the first assertion in four steps.

 \begin{stepp}
$k\btensor_RD$ is acyclic.
 \end{stepp}

Let $Y\to M$ be a flat resolution.  For $V=\hh{k\otimes_RY}$ one
has $V_0=M/\fm M\ne 0$, so $k$ is a direct summand of $V$, hence
it suffices to prove $V\btensor_RD$ is acyclic.

Let $X\to k$ be a flat resolution.  As $k\otimes_RY$ is a complex of $k$-vector 
spaces, one may choose the first one of the quasi-isomorphisms below:
 \[
 V \xra{\ \simeq\ }
 k\otimes_RY \xla{\ \simeq\ }
 X\otimes_RY\,.
 \]
The second one is standard. Proposition
\ref{flags:semi-flat}(1) and formula \eqref{btensor:associative}
now yield
 \begin{align*}
 \hh{V\btensor_RD}
 &\cong\hh{(k\otimes_RY)\btensor_RD}
\\
 &\cong\hh{(X\otimes_RY)\btensor_RD}
\\
 &=\hh{X\btensor_R(Y\btensor_RD)}\,.
 \end{align*}

Proposition \ref{flags:semi-flat}(1) also gives the first
quasi-isomorphism below:
 \[
 Y\btensor_RD\xra{\ \simeq\ }M\btensor_RD\xra{\ \simeq\ }0\,.
 \]
The second one is our hypothesis.  {}From Proposition
\ref{diff:quisms}(1) we get isomorphisms
 \[
 \hh{X\btensor_R(Y\btensor_RD)}\cong
 \hh{X\btensor_R(M\btensor_RD)}\cong
 \hh{X\btensor_R 0} =0\,.
 \]
 The two chains of isomorphisms yield $\hh{ V\btensor_RD}=0$, as desired.

 \begin{stepp}
$L\btensor_RD$ is acyclic for each $R$-module $L$ of finite length.
 \end{stepp}

We induce on $\length_RL$.  When it is $1$ one has $L\cong k$, so
the desired result was established in Step 1. For $\length_RL\ge2$ there
is an exact sequence of $R$-modules $0\to k\to L\to L'\to 0$. It
induces an exact sequence of differential $R$-modules
 \[
0\lra k\btensor_RD\lra L\btensor_RD\lra L'\btensor_RD\lra 0\,.
 \]
As $\length_R L'= \length_RL-1$, the induction hypothesis and the exact
triangle \eqref{cone:triangle} yield $\hh{L\btensor_RD}=0$. This
completes the proof of step 2.

 \begin{stepp}
$W\btensor_RD$ is acyclic for each bounded complex $W$ of $R$-modules,
such that the $R$-module $\HH hW$ has finite length for each $h$.
 \end{stepp}

Set $i=\inf\{h\var \HH hW\ne 0\}$.  The inclusion into $W$ of the
subcomplex
 \[
\quad\cdots \lra W_{i+2}\lra  W_{i+1}\lra \Ker(\dd_i)\lra 0\lra 0
\lra\cdots
 \]
is a quasi-isomorphism, so by Proposition \ref{flags:semi-flat}(1)
we may assume $W_h=0$ for $h<i$.

Set $H=\shift^i\HH iW$, let $\pi\col W\to H$ be the augmentation,
and let $j$ be the number of non-zero homology modules of $W$. The
exact sequence of complexes
 \[
0\lra\Ker(\pi)\lra W \lra H\lra 0
 \]
shows $\Ker(\pi)$ has $j-1$ non-vanishing homology modules of finite
length.  Since one has $\hh{H\btensor_RD}=0$ by Step 2, the 
exact sequence of differential modules
 \[
0\lra \Ker(\pi)\btensor_RD\lra W\btensor_RD\lra H\btensor_RD\lra0
 \]
yields a homology triangle \eqref{diff:triangle}, from which
$\hh{ W\btensor_RD}=0$ follows by induction.

 \begin{stepp}
$D$ is acyclic.
 \end{stepp}

Let $K$ be the Koszul complex on a finite generating set for the maximal
ideal of $R$.  Step 3 shows that $K\btensor_RD$ is acyclic, hence so
is $D$, by Lemma \ref{koszul} below.
 \end{proof}

\begin{lemma}
 \label{koszul}
Let $(R,\fm,k)$ be a local ring and $K$ the Koszul
complex on elements $x_1,\dots,x_e$ in $\fm$.  Let  $D$ be a
differential $R$-module with $\hh D$ finitely generated.

The module $\hh{K\btensor_RD}$  is then finitely generated.

Furthermore, one has $\hh{K\btensor_RD}=0$ if and only if $\hh D=0$.
 \end{lemma}

 \begin{proof}
Recall that $K$ is a tensor product $K'\otimes_RK''$, where $K'$
and $K''$ are Koszul complexes on the sequences $x_1$ and
$x_2,\dots,x_e$, respectively.  Thus, one has a canonical
isomorphism $K\btensor_RD=K'\btensor_R(K''\btensor_RD)$ of
differential modules, see \eqref{btensor:associative}.  By
induction, it suffices to prove the lemma for $e=1$, so we set
$x=x_1$.

The Koszul complex on $x$ is the cone of $x\idmap^R$. Using
\eqref{btensor:compcones} and \eqref{btensor:leftunit} one gets
 \[
\cone{x\idmap^R}\btensor_RD=\cone{(x\idmap^R)\btensor_RD}=
\cone{x\idmap^{D}}\,.
 \]
Thus, the exact triangle \eqref{cone:triangle} gives an exact
sequence of $R$-modules
 \[
0\lra\hh D/x\hh D\lra\hh{K\btensor_RD}\lra\hh D
 \]
The assertions of the lemma follow, with a nod from Nakayama for
the last one.
 \end{proof}

We describe a final tool needed for the proof of Theorem
\ref{intersection:local}.

 \begin{bfchunk}{{\v C}ech complexes.}
 \label{stables}
For each element $x\in R$, the localization map $R\to R_x$
defines a complex of $R$-modules with $R$ in degree $0$, as
follows:
 \[
C(x)=\quad \cdots\lra0\lra R\lra R_x\lra 0\lra\cdots
 \]
Let $\bsx=x_1,\dots,x_{s+1}$ be a sequence of elements in $R$. Set
 \begin{equation}
 \label{stables:complexes}
\cpxx Cn =  \begin{cases}
R & \text{if $n=0$}\\
C(x_1) \otimes_R \cdots \otimes_R C(x_n) & \text{if $n\geq 1$.}
 \end{cases}
 \end{equation}
The complex $C^{n}$ is concentrated between degrees $-n$ and $0$.
It is the modified {\v C}ech complex on $x_1,\dots,x_n$; see \cite[\S
3.5]{BH}.

Since $\cpxx C{n+1} = \cpxx Cn\otimes_RC(x_{n+1})$, the exact
sequence of complexes
 \[
0\lra\shift^{-1}R_{x_{n+1}}\lra C(x_{n+1})\lra R\lra 0
 \]
induces an exact sequence of bounded complexes of flat $R$-modules:
 \begin{equation}
\label{stable:ses} 0\lra\shift^{-1}\cpxx Cn\otimes_RR_{x}\lra
\cpxx C{n+1} \xra{\cpxx{\eps}{n+1}} \cpxx Cn \lra0\,.
 \end{equation}
 \end{bfchunk}

 \begin{proof}[Proof of Theorem \rm{\ref{intersection:local}}]
Set $s=\dim R-1$.  Let $M$ be a big Cohen-Macaulay $R$-module and
let $\bsx=x_1,\dots,x_{s+1}$ be a system of parameters that forms an
$M$-regular sequence, see \eqref{bigCM}.  Replacing the $x_i$
with their powers, we assume $\bsx\hh D=0$.

For the entire proof we fix the complexes $\cpxx
X{n}= M\otimes_R\cpxx C{n+1}$, obtained for $n=-1,\dots,s$ from
\eqref{stables:complexes}, and use \eqref{stable:ses} to define morphisms
 \[
\cpxx{\vartheta}{n}\col\cpxx X{n} \xra{\ M\otimes_R\cpxx \eps{n}\
}\cpxx X{n-1} \quad\text{for}\quad n=0,\dots,s\,.
 \]

The $M$-regularity of $\bsx$ implies $\HH i{\cpxx Xn}=0$ for $i\ne n$,
see \cite[(3.5.6) and (1.6.16)]{BH}.  Therefore, one has inclusions
\[
\Image(\HH i{\cpxx{\vartheta}{n}})\subseteq\HH i{\cpxx X{n-1}}=0\,.
\]
They yield equalities
 \begin{equation}
 \tag{a}
\hh{\cpxx{\vartheta}{n}}=0  \quad\text{for}\quad n=0,\dots,s\,.
 \end{equation}

By hypothesis, there is a split epimorphism $\pi\col F\to D$ of
differential $R$-modules, where $\hh D$ has finite non-zero
length. Proposition \ref{carlsson} shows that the condition
 \begin{equation}
 \tag{b}
\hh{\vartheta\btensor_R\pi}\ne0 \quad\text{for}\quad
\cpxx{\vartheta}{0}\circ\cdots\circ\cpxx{\vartheta}{s}\,.
  \end{equation}
gives the desired conclusion: $\flclass RF\ge s+1$.  We verify (b)
in three steps.

\setcounter{stepp}{0}

 \begin{stepp}
$(\cpxx X{n}\otimes_RR_{x_n})\btensor_R D$ is acyclic for
$n=-1,\dots,{s}$.
 \end{stepp}

Fix $n$, set $x=x_{n}$, and choose a flat resolution $Y\to M$. The
arrow below
 \[
Y\otimes_R\cpxx C{n+1}\otimes_RR_x \lra M\otimes_R \cpxx C{n+1}\otimes_RR_x
=\cpxx X{n}\otimes_RR_x
 \]
is a quasi-isomorphism because $\cpxx C{n+1}\otimes_RR_x$ is a
bounded complex of flat modules. By Proposition
\ref{flags:semi-flat}(1), it induces the first quasi-isomorphism
of differential modules
 \begin{align*}
  (\cpxx X{n}\otimes_RR_{x})\btensor_R D &\simeq (Y\otimes_R \cpxx
  C{n+1}\otimes_RR_{x})\btensor_R D\\ &\cong  (Y\otimes_R \cpxx
  C{n+1})\btensor_R(R_{x}\otimes_R D)\,.
 \end{align*}
The isomorphism is due to associativity.  As $R_x$ is $R$-flat and
$x\cdot\hh D=0$, we get $\hh{R_{x}\otimes_RD}\cong R_{x}\otimes_R\hh{D}=0$.
That is, the map $R_x\otimes_RD\to0$ is a quasi-isomorphism.  As
$Y\otimes_R \cpxx C{n+1}$ is a bounded below complex of flat
modules, it induces
 \[
  (Y\otimes_R \cpxx C{n+1})\btensor_R(R_{x}\otimes_R D) \simeq  (
  Y\otimes_R \cpxx C{n+1})\btensor_R0 =0\,,
 \]
see Proposition \ref{diff:quisms}.  The proof of Step 1 is now complete.

 \begin{stepp}
$\hh{\cpxx{\vartheta}{n}\btensor_RD}\col\hh{\cpxx
X{n}\btensor_RD}\to \hh{\cpxx X{n-1}\btensor_RD}$ is bijective for
$n=0,\dots,s$.
  \end{stepp}

The complexes in the exact sequence \eqref{stable:ses} consist of
flat $R$-modules, so
 \[
0\lra\shift^{-1}(\cpxx X{n-1}\otimes_RR_{x})\lra
 \cpxx X{n}\xra{\ \cpxx{\vartheta}{n}\ }\cpxx X{n-1}\lra0
 \]
is an exact sequence of complexes of $R$-modules.  The induced
sequence
 \[
0\to\shift^{-1}(\cpxx X{n-1}\otimes_RR_{x})\btensor_R D\to
  \cpxx X{n}\btensor_RD \xra{\ \cpxx{\vartheta}{n}\btensor_RD\ }
   \cpxx X{n-1}\btensor_RD\to 0
 \]
of differential $R$-modules is exact. Its homology exact triangle,
\eqref{diff:triangle} and the result of Step 1 imply that
$\hh{\cpxx{\vartheta}{n}\btensor_RD}$ is bijective. Step 2 is
complete.

 \begin{stepp}
Condition (b) holds.
  \end{stepp}

Indeed, the map $\hh{\vartheta\btensor_R\pi}$ factors as a
composition
 \[
\xymatrixrowsep{2pc} \xymatrixcolsep{5pc} \xymatrix{ \hh{\cpxx
X{s}\btensor_RF}\ar@{->}[r]^-{\hh{\cpxx X{s}\btensor_R\pi}}
  &\hh{\cpxx X{s}\btensor_RD}\ar@{->}[r]^-{\hh{\vartheta\btensor_RD}}
  &\hh{\cpxx X{-1}\btensor_RD}}\,.
 \]
The first map above is surjective because $\pi$ is a split
epimorphism of differential modules.  The second map is a
composition of isomorphisms, due to the result of Step 2.
Finally, the module $\hh{\cpxx X{-1}\btensor_RD}$ is not zero, 
by Theorem \ref{trivialflags:local}.
 \end{proof}

\section{Class inequality. II}
\label{Class inequality. II}

In this section we prove global, relative versions of the results in
Section \ref{Class inequality. I}.
The theorem below can be compared with Theorem \ref{intersection:local},
through Remark \ref{class}(8).

 \begin{theorem}
 \label{intersection}
Let $R\to S$ be a homomorphism of commutative
noetherian rings.  Let $F$ be a finitely generated differential
$R$-module and $D$ a retract of $F$.  Let $\fq$ be a prime ideal
in $S$ minimal over $IS$, where $I=\Ann_R\hh D$.

When $S_\fq$ has a big Cohen-Macaulay module one has
\[
  \prclass RF\geq \dim S_\fq\,.
\]

An inequality $\prclass RF\geq \dim S_\fq -1$ holds in general.
 \end{theorem}

 Recall a notion introduced by Hochster \cite{Ho:LNM}: for an ideal $I$ in $R$, set
 \begin{gather*}
  \supheight I = \sup\left\{\height(IS)\left|
   \begin{gathered}
 \text{$R\to S$ is a homomorphism}\\
 \text{of rings and $S$ is noetherian}
   \end{gathered}
\right\}\right.\,.
 \end{gather*}
Evidently, one has $\supheight I\geq \height I$, whence the notation.

Every local ring containing a field has a big Cohen-Macaulay module,
cf.\ \eqref{bigCM}, so the next result is essentially a reformulation
of part of the theorem.  The second inequality in it implies the Class
Inequality, stated in the introduction.

\begin{corollary}
\label{intersection:global}
 \pushQED{\qed}%
If $R$ is a commutative noetherian ring, $F$ a finitely
generated differential $R$-module, $D$ a retract of $F$,
and $I=\Ann_R \hh D$, then one has
   \[
\prclass RF \geq \supheight I-1\,.
   \]

When $R$ is an algebra over a field a stronger inequality holds:
   \[
\prclass RF \geq \supheight I\,.
   \]

If $\dim R\leq 3$ holds, or if $R$ is Cohen-Macaulay, then one has
   \[
\prclass RF \geq \height I\,.
 \qedhere
   \]
 \end{corollary}

Hochster's motivation for introducing super heights was to prove the
following homological generalization of Krull's Principal Ideal Theorem:
If $R$ is a noetherian ring containing a field, then every finitely
generated $R$-module $M$ satisfies 
\[
\supheight(\Ann_RM)\le\pd_RM\,;
\]
see \cite{Ho:LNM}, also \cite[(9.4.4)]{BH}.  In view of Remark \ref{class}(6), this result
may be recovered by applying the corollary to $F=\dfm P$, where $P$ is a projective
resolution of $M$, of length equal to $\pd_RM$.

Examples show that the last two inequalities in Corollary
\ref{intersection:global} are sharp:

\begin{example}
\label{optimal}
Let $R$ be a commutative noetherian algebra over a field. For each integer
$d$ with $0\leq d\leq \dim R$, there exists a differential $R$-module $F$
for which
 \begin{gather*}
  \prclass RF =d = \supheight I \quad\text{where $I=\Ann_R \hh F$.}
 \end{gather*}

 Indeed, fix such a $d$. One can then find an ideal $I$ of height $d$ generated by a set
 $\bsx$ with $d$ elements.  Let $K$ be the Koszul complex on $\bsx$, and set $F=\dfm
 K$. Since $I\HH iK=0$ for each $i$ and $\HH 0K=R/I$, one has $\Ann_R\hh F=I$. Remark
 \ref{class}(6) implies the first inequality below:
 \[ 
d\geq \prclass RF\geq \supheight I \geq \height I = d\,.
 \]
The second one is given by Corollary \ref{intersection:global};
 the last one holds  always.
 \end{example}

To prove Theorem \ref{intersection} we need information on how the
support of the homology of differential modules is affected by
base change.  The next result provides a complete and completely
satisfying answer in the presence of projective flags.

 \begin{theorem}
 \label{flags:basechange}
Let $R\to S$ be a homomorphism of noetherian commutative rings and
$D$ a differential $R$-module with $\hh D$ finitely generated.

If $D$ is a retract of some differential module admitting a projective
flag, then
 \[
\Supp_S\hh{S\otimes_RD}=\Supp_S(S/IS) \quad\text{where}\quad
I=\Ann_R\hh D\,.
 \]
 
When, in addition, the $S$-module $S\otimes_RD$ is finitely
generated, the homology module $\hh{S\otimes_RD}$ has finite
length over $S$ if and only if the ring $S/IS$ is artinian.
\end{theorem}

\begin{proof}
Let $\fq$ be a prime ideal of $S$, set $\fp=R\cap \fq$, and
let $R_\fp\to S_\fq$ be the induced local homomorphism.  One then
has isomorphisms 
\[
\hh{S\otimes_RD}_\fq \cong \hh{(S\otimes_RD)_\fq} \cong
 \hh{S_{\fq}\otimes_{R_\fp}D_\fq}
\]
of $S_\fq$-modules.  Since the differential $R_\fp$-module $D_\fp$ and the $R_\fp$-module
$S_\fq$ satisfy the hypotheses of Theorem \ref{trivialflags:local}; it shows that
$\hh{S_{\fq}\otimes_{R_\fp}D_\fq}=0$ holds if and only if $\hh{D_\fp}=0$ does, that is to
say, if and only if $\hh{D}{}_\fp=0$. As $\hh D$ is finite over $R$, the last condition is
equivalent to $\fp\nsupseteq I$, which is tantamount to $\fq\nsupseteq IS$.  We have now
established that the $S$-modules $\hh{S\otimes_RD}$ and $S/IS$, have the same support.

If $S\otimes_RD$ is finitely generated over $S$, then so is
$\hh{S\otimes_RD}$, hence its length is finite if and only if
its support consists of maximal ideals.  This support being equal to
that of $S/IS$, the last condition is equivalent
to the ring $S/IS$ being artinian. \end{proof}

 \begin{proof}[Proof of Theorem \rm{\ref{intersection}}]
 \label{proof:intersection}
When $S_\fq$ has a big Cohen-Macaulay module, one gets
 \[
\prclass RF \geq \prclass{S_\fq}{(S_\fq\otimes_RF)}
 \]
from Remark \ref{class}(10).  The differential $S_\fq$-module
$S_\fq\otimes_RD$ is a retract of $S_\fq\otimes_RF$.  Since the
length of $S_\fq/IS_\fq$ is non-zero and finite, Theorem
\ref{flags:basechange} implies that the same holds for the length
of $\hh{S_\fq\otimes_RD}$. Therefore, Theorem
\ref{intersection:local} yields
 \[
\prclass{S_\fq}{(S_\fq\otimes_RF)}\ge\dim S_\fq\,.
 \]
The preceding inequalities yield $\prclass RF\geq
\dim S_\fq$, as needed.

Next we drop the assumption on $\fq$.  Let $p$ denote the characteristic of the residue
field of $S_\fq$.  The ring $S'=S/pS$ is an algebra over the prime field of characteristic
$p$, so it has a big Cohen-Macaulay module, see Remark \ref{bigCM}.  The already
established part of the theorem yields the first inequality:
\[
\prclass{R}{F}\ge\dim S'\geq\dim S_\fq-1\,.
\]
The second one holds always. This completes the proof.
 \end{proof}

\section{Rank inequalities}
\label{Rank inequalities}

The results of Sections \ref{Class inequality. I} and \ref{Class inequality. II} provide
lower bounds on the projective class of a differential module in terms of invariants of
its homology. Here we provide similar bounds for the rank of a differential module
admitting a finite free flag.

\begin{theorem}
\label{horrocks:minor}
Let $R$ be a commutative noetherian algebra over a field.

If $D$ is a retract of a differential $R$-module $F$ that admits a finite
free flag, then
 \[
\rank RF\geq 2(\supheight I) \quad\text{where}\quad I=\Ann_R \hh D\,.
 \]
\end{theorem}

\begin{Remark}
Not surprisingly, the proof shows that over every noetherian ring
one has an inequality $\rank RF\geq 2(\supheight I)-2$.
\end{Remark}

When the height of $\Ann_R\hh F$ is at least $5$, this remark implies $\rank RF\geq
8$. For smaller heights we have the following result.

\begin{theorem}
\label{horrocks:ufd}
Let $R$ be a commutative noetherian ring, $F$ a differential $R$-module, and set
$I=\Ann_R\hh F$ and $d=\height I$.

If $F$ admits a finite free flag, then one has inequalities:
 \begin{alignat}{2}
 \label{ufd:small}
\rank RF &\geq 2d&\quad\text{when}\quad d&\leq 3\,;
 \\
 \label{ufd:large}
\rank RF &\geq 8&\quad\text{when}\quad d&\geq 3
\text{ and $R$ is a UFD.}
 \end{alignat}
\end{theorem}

For differential graded modules with finite length homology over graded polynomial rings,
see Remark \ref{CT} for details, our theorem specializes to \cite[Thm.\ 2]{Ca:ams}.
Theorems \ref{horrocks:minor} and \ref{horrocks:ufd} are proved at the end of the section.
Together, they contain the Rank Inequality, stated in the introduction, and suggest
the following:

\begin{conjecture}
 \label{horrocks:major}
Let $R$ be a local ring and $F$ a differential $R$-module with a finite
free flag. If $\hh F$ has non-zero  finite length, then 
\[
\rank RF\geq 2^d\quad \text{for} \quad d={\dim R}\,.
\]
 \end{conjecture}

This is in line with several results and open problems in algebra and topology:

\begin{remark}
\label{horrocks:beh}
Buchsbaum and Eisenbud \cite[(1.4)]{BE}, and Horrocks \cite[Pbl.\ 24]{Ha} conjectured that
if $P$ is a finite free resolution of a module of finite length over a local ring $R$ of
dimension $d$, then $\rank R{P_n}\ge \binom {d}n$. These inequalities predict
 \[
\rank R(\dfm P)=\sum_n \rank R{P_n} \geq 2^d\,,
 \]
as does Conjecture \ref{horrocks:major}.  For $d\le 4$ the
Buchsbaum-Eisenbud-Horrocks conjecture follows from the Generalized
Principal Ideal Theorem. For equicharacteristic rings
Avramov and Buchweitz \cite[(1)]{AB} use Evans and Griffith's Syzygy
Theorem \cite{EG} to prove $\rank R(\dfm P)\geq\frac32(d-1)^2 + 8$
for $d\ge5$; thus, for $d=5$ one has $\rank R(\dfm P)\geq2^5$.
 \end{remark}

 \begin{remark}
 \label{horrocks:cap}
Let $X$ be a finite CW complex.  Halperin \cite[(1.4)]{Hl} asked: If
the torus $(\BR/\BZ)^d$ acts with finite isotropy groups on $X$, does
then one have
\[
\sum_n \rank{\BQ}{\HH n{X;\BQ}}\ge 2^d\,?
\]
In a similar vein, Carlsson \cite[Conj.\ I.3]{Ca:ln} conjectured that if the
elementary abelian group $(\BZ/2\BZ)^d$ acts freely on $X$, then one has
\[
\sum_n \rank{\BF_2}{\HH n{X;\BF_2}}\ge 2^d\,.
\]
It is known that if Conjecture \ref{horrocks:major} holds for differential
graded modules over a graded polynomial ring in $d$ variables over
a field, then so do the inequalities above; see the remarks after
\cite[Conj.\ II.2]{Ca:ln}.  In particular, Theorem \ref{horrocks:ufd}
implies certain results of Allday and Puppe; see \cite[(1.4.21) and
(4.4.3)(1)]{AP}.
 \end{remark}

We note that the conclusion of Theorems \ref{horrocks:ufd} and
\ref{horrocks:minor} may fail if the hypothesis on $F$ is weakened from
admitting a free flag to just being free as an $R$-module:

 \begin{example}
 \label{notaflag}
Let $R=k[x,y]$ be a polynomial ring over a field $k$ and $D$ the differential $R$
module given by the square-zero matrix
 \[
A=\begin{bmatrix} xy & -x^2\\ y^2 & -xy\end{bmatrix}
 \]
  A straightforward
calculation yields
\begin{gather*}
 \Ker(\delta) = R\begin{bmatrix} x \\ y
\end{bmatrix}
 \quad\text{and}\quad
 \Image(\delta)
 =(x,y)R \begin{bmatrix}
x
 \\
y
 \end{bmatrix}
\end{gather*}
Therefore, one has $I = (x,y)$, as well as $\rank RF=2<4=2\height I$.
\end{example}

The proofs of the preceding theorems require preparation.

When $R$ is a domain, $R_0$ its field of fractions, and $M$ an $R$-module
one sets $\rank RM=\rank {R_0}({R_0}\otimes_RM)$.  For a matrix $A$ with
entries in $R$ let $\rank R(A)$ denote its rank over $R_0$.  Rank is used
to form Euler-Poincar\'e characteristics of complexes.  Only vestigial
versions of Euler's formula hold in the absence of gradings.

\begin{remark}
If $R$ is an integral domain and $D$ is a finitely generated
differential $R$-module, then the following hold:
 \begin{align}
  \label{rank:double}
 \rank RD&=2\rank R(\delta^D)\iff
 \rank R\hh{D}=0
 \\
  \label{rank:parity}
 \rank RD&\equiv \rank R{\hh{D}}\pmod 2\,.
\end{align}

Indeed, both formulas result from the equality
 \[
\rank RD=\rank R{\hh D}+2\rank R{\Image(\delta^D)}
 \]
obtained by additivity of rank from the exact sequences
of $R$-modules
 \[
\xymatrixrowsep{.3pc} \xymatrixcolsep{1pc} \xymatrix{
 0\ar@{->}[r]&\Ker(\delta^D)\ar@{->}[r]&D\ar@{->}[r]&\Image(\delta^D)\ar@{->}[r]&0
\\
 0\ar@{->}[r]&\Image(\delta^D)\ar@{->}[r]&\Ker(\delta^D)\ar@{->}[r]&\hh D\ar@{->}[r]&0\,.
 }
 \]
\end{remark}

\begin{lemma}
\label{rank:domain}
If $R$ is a domain and $F$ is a finitely generated differential
$R$-module with $\frclass RF=l<\infty$, then one has
 \begin{equation}
  \label{rank:bound}
 \rank RF\geq 2l\,.
 \end{equation}
Moreover, if $\Filt F{n}$ is a free flag  with $F^l=F$, then
the following hold:
 \begin{align}
 \label{rank:layer1}
 \rank R(F_n) &\ge
 \begin{cases}
1 & \text{for $n=0$ or $n=l$} \cr 2 &\text{for  $n=1,\dots, l-1$}
 \end{cases}
\\
 \label{rank:layer2}
 \delta(\filt Fn)&\not\subseteq\filt F{n-2}
 \quad\text{for}\quad n=1,\dots, l\,.
 \end{align}
 \end{lemma}

 \begin{proof}
Assuming \eqref{rank:layer2} fails for some $n$, one finds in $F$
a new flag $\Filt Gi$ by setting
 \begin{gather*}
\filt Gi = \begin{cases}
\filt F{i}  & \text{for $i\leq n-2$} \\
\filt F{i+1} & \text{for $i\ge n-1$}
 \end{cases}
 \end{gather*}
It implies $\frclass RF<l$, contradicting our hypothesis.  Thus,
\eqref{rank:layer2} holds.

Assume next that \eqref{rank:layer1} fails for some $n$.  We show
that \eqref{rank:layer2} fails for some $j$, and thus draw a
contradiction.  For $n=0$, $n=l$, or $F_n=0$, one may take $j=n+1$.  For
$n\in[1,l-1]$ the complex \eqref{spectral sequence} has the form
 \[
\xymatrixcolsep{1.5pc} \xymatrix{ \cdots \ar@{->}[r]  &
F_{n+1}\ar@{->}[rr]^-{\dd_{n,n+1}}
    &&R  \ar@{->}[rr]^-{\dd_{n-1,n}} && F_{n-1}\ar@{->}[r] &\cdots}
 \]
Since $R$ is a domain, either $\dd_{n+1}=0$ or $\dd_n=0$; set
$j=n+1$ or $j=n$, respectively.

Finally, formula \eqref{rank:bound} follows from the computation
 \begin{align*}
\rank RF &= \rank R{(F_0)} + \sum_{n=1}^{l-1} \rank R{(F_n)} +
\rank R{(F_l)}\\
         & \geq 1 + (l-1)2 + 1\\
         & =2l\,,
 \end{align*}
where the inequality comes from the already established formula
\eqref{rank:layer1}.
 \end{proof}

 \begin{proof}[Proof of Theorem \rm{\ref{horrocks:minor}}]
 \label{proof:minor}
Let $R\to S$ be a homomorphism to a noetherian ring $S$, such that
$\supheight I=\height(IS)$.  Pick a prime ideal $\fq$ in $S$,
minimal over $IS$, then a prime ideal $\fp\subseteq\fq$, such that
$\height(IS)=\dim(S_\fq)=\dim(S_\fq/\fp S_\fq)$.

Let $R'$ denote the local domain $S_\fq/\fp S_\fq$.  The
differential graded $R'$-module $D'=R'\otimes_RD$ is a retract of
$F'=R'\otimes_RF$, so the following inequalities
 \[
\rank {R'}{F'}\geq2\frclass{R'}{F'}\geq2\prclass{R'}{F'}
\geq2\dim{R'}
 \]
are given by formula \eqref{rank:bound}, Remark \ref{class}(10),
and Theorem \ref{intersection:local}.   The desired result
follows, as one has $\rank RF=\rank {R'}{F'}$ and $\dim
R'=\supheight I$.
 \end{proof}


 \begin{proof}[Proof of Theorem \rm{\ref{horrocks:ufd}}]
Recall that $F$ is a differential $R$-module admitting a finite
free flag, $\Ann_R\hh F=I\ne R$, and $d=\height I$.

When $d=1$ the desired  inequality $\rank RF\ge 2d$ clearly
holds.

When $d=2$ or $d=3$ pick prime ideals $\fq\supseteq I$ and $\fp\subseteq \fq$ such that
$\height (\fq/\fp)=2$, and set $S= R_\fq/\fp R_\fq$. Note that $\length_S\hh{S\otimes_RF}$
is finite and non-zero, by Theorem \ref{flags:basechange}, and that $S$ has a big
Cohen-Macaulay module, since $\dim S\leq 3$.  Applying successively Remark \ref{class}(9),
formula \eqref{rank:bound}, and Theorem \ref{intersection} we obtain
\[
\rank RF\geq \rank S{(S\otimes_RF)} \geq 2\prclass
S{(S\otimes_RF)} \geq 2\cdot d\,.
\]
This completes the proof of \eqref{ufd:small}.

Assume $d\geq 3$ and $R$ is a UFD.  Let $\fq$ be a prime ideal
containing $I$ such that $\dim R_\fq=d$.  Since $R$ is a domain,
\eqref{rank:bound} yields the first inequality below
\[
\rank RF \geq 2\prclass RF \geq 2\cdot\min\{3,d-1\} = 6\,.
\]
The second is due to Theorem \ref{intersection}, because when
$d=3$ the ring $R_\fq$ has a big Cohen-Macaulay module, see
\ref{bigCM}.  Formula \eqref{rank:parity} rules out $\rank RF=7$,
so to finish the proof of the inequality $\rank RF\ge 8$ it
remains to show $\rank RF\ne 6$.

Replacing $R$ with $R_\fq$, for the rest of the proof we assume that $R$
is a local UFD, $\dim R\ge 3$, and $\hh F$ has finite non-zero length.
We assume $\rank RF= 6$ and draw a contradiction. Set $l=\frclass RF$
and let $\Filt Fn$ be a free flag with $F^l=F$. Lemma \ref{rank:domain}
implies that the first page of the spectral sequence in Remark
\ref{spectral sequence} is
 \[
\Gr 01{}Fn=\xymatrixcolsep{1.5pc} \xymatrix{
 \cdots\ar@{->}[r] &0\ar@{->}[r] & R
\ar@{->}[r]^{\dd_{23}} &R^2\ar@{->}[r]^{\dd_{12}}
    &  R^2 \ar@{->}[r]^{\dd_{01}} & R \ar@{->}[r] & 0\ar@{->}[r] &\cdots}
\]
with $\dd_{n,n+1}\ne0$ for $n=0,1,2$.  Let $A_n$ denote the matrix
of $\dd_{n,n+1}$ in some bases; $\dd_{23}\ne0$ implies $\rank R
{A_1}=1$.  Let $J$ be the ideal generated by the entries of $A_1$.

If $J=R$, then (possibly after changing bases) one can find
$x,z\ne0$ such that
 \[
 A_2=\begin{bmatrix} 0\\z\end{bmatrix}\,,\quad
 A_1=\begin{bmatrix} 1 &0\\0& 0\end{bmatrix}\,,\quad
 A_0=\begin{bmatrix} 0 & x\end{bmatrix}
 \]
 It follows that on the second page the only non-trivial complex is
 \[
\Gr 02{}Fn=\quad
\xymatrixcolsep{1.5pc} \xymatrix{
 \cdots\ar@{->}[r] &0\ar@{->}[r] & R/(z)
\ar@{->}[r]^{\dd_{02}} &R/(x)\ar@{->}[r] & 0\ar@{->}[r] &\cdots}
 \]
where $\dd_{02}$ is induced by multiplication with some $y\in(x:z)$.
The condition $l=3$ implies $\Gr{}3{}Fn=\gr{\infty}Fn$, see
\eqref{spectral:stablepage}, so \eqref{spectral:convergence} yields an
exact sequence
 \[
0\lra R/(x,y)\lra\hh F\lra (x:y)/(z)\lra0
 \]
of $R$-modules.  Since $R$ is a UFD, the ideal $(x:y)$ is
generated by $x'=x/v$ where $v=\operatorname{gcd}(x,y)$.  Thus,
the module $(x:y)/(z)$ is isomorphic to $R/(w)$, where $w=z/x'$.
The hypothesis $\hh F\ne0$ implies $R/(w)\ne0$ or $R/(x,y)\ne0$.
In the first case we get $\dim_RR/(w)= \dim R-1\ge2$, in the
second $\dim_RR/(x,y)\ge\dim R-2\ge1$.  Either inequality
contradicts the hypothesis that $\hh F$ has finite length.

If $J\ne R$, choose a prime ideal $\fp$ minimal over $J$, and let $S$
denote the field $R_\fp/\fp R_\fp$. Corollary \ref{heights} implies
$\height\fp\leq 2$, so one has $\fp\nsupseteq I$, and thus $\hh{F_\fp}=0$.
Theorem \ref{trivialflags:local} now yields $\hh{S\otimes_RF}=0$.
The first page of the spectral sequence $\eGr pr{}=\Gr pr{}{S\otimes_RF}n$
associated to the flag $\Filt {S\otimes_RF}n$ is the complex
 \[
\eGr 01{}=\quad\xymatrixcolsep{1.7pc} \xymatrix{
 \cdots\ar@{->}[r] &0\ar@{->}[r] & S
\ar@{->}[r]^{\dd_{23}} &S^2\ar@{->}[r]^{0}
    & S^2 \ar@{->}[r]^{\dd_{01}} & S \ar@{->}[r] & 0\ar@{->}[r] &\cdots}
\]
of vector spaces over $S$.  The two complexes on the second page
$\eGr {}2{}$ are
 \[
  \xymatrixrowsep{.3pc} \xymatrixcolsep{1.8pc} \xymatrix{
  \eGr 02{}=
\\
\eGr 12{}=
}
 \quad
\xymatrixrowsep{.3pc} \xymatrixcolsep{1.8pc} \xymatrix{
   \cdots\ar@{->}[r]
  &0\ar@{->}[r]
  &\Coker(\dd_{23})\ar@{->}[r]^-{\dd_{02}}
  &\Coker(\dd_{01})\ar@{->}[r]
  &0\ar@{->}[r]
  &\cdots
\\
   \cdots\ar@{->}[r]
  &0\ar@{->}[r]
  &\Ker(\dd_{23})\ar@{->}[r]^-{\dd_{13}}
  &\Ker(\dd_{01})\ar@{->}[r]
  &0\ar@{->}[r]
  &\cdots
 }
\]
Counting ranks over $S$, one verifies the following assertions:
$\eGr 232\ne 0$ when $\dd_{23}=0$, or when $\dd_{23}\ne0$ and
$\dd_{01}\ne0$ hold simultaneously; $\eGr 131\ne 0$ when
$\dd_{23}\ne0$, but $\dd_{01}=0$.  Now \eqref{spectral:stablepage}
yields $\eGr 2{\infty}2 = \eGr 232$ and $\eGr 1{\infty}1 = \eGr
131$, so one has $\hh{S\otimes_RF}\ne0$ by
\eqref{spectral:convergence}. This new contradiction completes the
proof of \eqref{ufd:large}.
 \end{proof}

 \section{Square-zero matrices}
 \label{Square-zero matrices}

Our primary purpose in this short section is to record matrix versions of
theorems proved earlier in the text.  As a side benefit we get a framework
for describing examples.  We assume that $R$ is an IBN (= invariant basis
number) associative ring: in every finitely generated free $R$-module
any two bases have the same number of elements, see \cite{Co-book}.
Commutative rings and left noetherian rings have this property.

Let $\operatorname{M}_s(R\op)$ be the ring of $s\times s$ matrices
over $R\op$.  Let $E$ be the $R$-module of $s\times1$ matrices with
entries in $R$ and $e_i\in E$ the matrix with $1$ in the $i$th row
and $0$ elsewhere.  For each $A=(a_{ij})\in\operatorname{M}_s(R\op)$,
define $\eps_A\in\End R{E}$ by the condition
 \[
\eps_A(e_j)=\sum_{i=1}^s a_{ij}e_i
 \quad\text{for}\quad j=1,\dots,s
  \]
The map $A\mapsto\eps_A$ is an isomorphism of rings
$\operatorname{M}_s(R\op)\cong\End R{E}$. The standard operations of
$\operatorname{M}_s(R\op)$ turns $E$ into a bimodule, so  $\Ker(A)=\{e\in
E\var Ae=0\}$ is an $R$-submodule.  Let $\Image(A)$ be the $R$-submodule
of $E$ spanned by the columns of $A$.

The map $A\mapsto(E,\eps_A)$ induces a bijection between conjugacy
classes of square-zero matrices in $\operatorname{M}_s(R\op)$ and
isomorphism classes of free differential $R$-modules of rank $s$; one
has $\hh{E,\eps_A}\cong\Ker(A)/\Image(A)$.  A matrix $A$ with $A^2=0$
is conjugated to a strictly upper triangular matrix if and only if
$(E,\eps_A)$ is a free flag.

A ring is \emph{projective-free} if it has the IBN property and its
finitely generated projectives are free, see \cite{Co-book}.  The remarks
above translate Theorem \ref{flags:triviality} into the next statement,
where $0_r$ and $1_r$ denote the $r\times r$ zero and identity matrices,
respectively.

 \begin{theorem}
 \label{Flags:split}
Let $R$ be a projective-free ring and let $A=(a_{ij})$ be an $s\times s$
matrix with entries in $R$.  The following conditions are equivalent.
 \begin{enumerate}[\rm(i)]
 \item
$s=2r$ and $A$ is conjugated to the matrix
 \(
\begin{bmatrix}0_r & 1_r \\ 0_r & 0_r\end{bmatrix}
 \).
 \item
$A$ is conjugated to some strictly upper triangular matrix, and
$(x_1,\dots,x_s)\in R^s$ is a solution of the system of equations
 \[
\sum_{j=1}^sa_{ij}x_j=0\quad\text{for}\quad i=1,\dots,s
 \]
if and only if $x_j=\sum_{i=1}^sa_{ji}c_i$ for fixed $c_1,\dots,c_s\in R$
and $j=1,\dots,s$. \qed
 \end{enumerate}
  \end{theorem}

Let $B$ be an $s\times s$ strictly upper triangular matrix with
entries in $R$.  A \emph{block partition} of $B$ in $l$
\emph{steps} is a sequence of integers $1\le s_0<\cdots<s_l=s$, such
that
 \[
(b_{uv})_{\substack{s_i\le u<s_{i+1}\cr s_j\le v<s_{j+1}}}=0
\quad\text{for all}\quad i\ge j\,.
 \]
When such a partition exists one has 
\[
l\ge\frclass
R(E,\eps_B)\ge\prclass R(E,\eps_B)\,.
\]
Thus, Theorems \ref{intersection}, \ref{horrocks:minor}, and
\ref{horrocks:ufd} translate into:

 \begin{theorem}
 \label{matrix:size}
Let $R$ be a commutative noetherian algebra containing a field and
let $A$ be an $s\times s$ matrix with entries in $R$, such that
$A^2=0$.  Let $I$ denote the ideal $\Ann_R\hh A$, assume $I\ne R$,
and set $d=\height I$.

When $A$ is conjugated to a strictly upper triangular matrix $B$,
every block partition of $B$ has at least $d$ steps, and
the inequality $s\ge 2d$ holds.

When, in addition, $R$ is a UFD and $d\ge3$ holds, one has $s\ge8$.
\qed
\end{theorem}

\appendix

\section{Ranks of matrices}
\label{Inner rank}

In this appendix $R$ is a commutative ring, $A$ is an $m\times n$ matrix
with entries in $R$, and $I_r(A)$ denotes the ideal in $R$ generated by
all $r\times r$ minors of $A$.  As usual, we set $I_0(A)=R$ and $I_r(A)=0$
for $r>\min\{m,n\}$, so that $I_r(A)\supseteq I_{r+1}(A)$ holds for all
$r\ge0$. We compare two notions of rank for $A$.

 \begin{remark}
  \label{detvsinn}
The \emph{determinantal rank} of $A$, denoted $\detrank RA$, is the
largest integer $r\ge0$ such that $I_r(A)\ne0$.  The \emph{inner rank}
of $A$, denoted $\innrank RA$, is the least integer $s\ge0$ such that
$A$ can be written as a product $A=A'A''$ with an $n\times s$ matrix $A'$
and an $s\times m$ matrix $A''$, see \cite{Co-book}.  Standard linear
algebra yields
 \[
\detrank RA \leq \innrank RA
 \]
and shows that equality holds when $R$ is a field.  It follows that
when $R$ is an integral domain the ranks coincide if and only if
$\innrank RA=\innrank KA$, where $K$ is the field of fractions of $R$.
Domains over which this holds for all $A$ have been described
by multiple conditions; in particular, by the property that the kernel
of every homomorphisms $R^n\to R^m$ is a union of free submodules,
see \cite[(5.5.9)]{Co-book}.
 \end{remark}

The curious result below complements the criterion for agreement of ranks.

\begin{proposition}
\label{factorization}
An integral domain $R$ is factorial if and only if every non-empty set
of principal ideals of $R$ has a maximal element, and each matrix $A$
over $R$ with $\detrank RA\le 1$ satisfies $\detrank RA =\innrank RA$.
 \end{proposition}

{}From Krull's Principal Ideal Theorem one obtains an immediate
consequence:

 \begin{corollary}
 \label{heights}
 \pushQED{\qed}%
If $R$ is a noetherian UFD and $\detrank RA=1$, then either
 \begin{equation*}
{\phantom{\square}} I_1(A)=R \quad\text{or}\quad
\height I_1(A) \leq \max\{m,n\}\,. \qedhere
 \end{equation*}
 \end{corollary}

 \begin{proof}[Proof of Proposition~\emph{\ref{factorization}}]
Let first $R$ be a factorial domain.  The maximality
condition on principal ideals holds by \cite[(0.9.4)]{Co-book}, so we
let $A$ be an $m\times n$ matrix over $R$ with $\detrank RA\le1$ and
show that its determinantal and inner ranks coincide.

There is nothing to prove if $\detrank RA=0$ or if ether $m$ or $n$
is equal to $1$, so we assume $A\ne0$ and $\min\{m,n\}\ge2$.  Let $B$
denote the $(m-1)\times n$ matrix consisting of the first $m-1$ rows
in $A$.  If $B=0$, then $A=A'A''$ holds for the $m\times1$ matrix $A'$,
transpose to $[0\ \cdots\ 0\ 1]$, and the $1\times n$ matrix $A''$ with
$a_j=a_{mj}$.  We further assume $B\ne 0$ and induce on $m$. If $m=2$,
then the matrix has the form
 \[
A=\begin{bmatrix}
x_{1} & \cdots & x_n \\
y_{1} & \cdots & y_n
\end{bmatrix}
 \]
with some $x_h\ne 0$.  As $\detrank R{(A)}\leq 1$, the rows of $A$
are proportional over the field of fractions $R_0$; in other
words, $y_j= (p/q)x_j$ for some $p/q\in R_0$ and $1\leq j\leq n$.
As $R$ is a UFD, one may assume $p,q$ are relatively prime. One
can then find $y_1',\dots,y_n'\in R$ satisfying $y_j=py_j'$ for
each $j$. This implies $x_j=qy'_j$, hence
\[
A= \begin{bmatrix} q\\p  \end{bmatrix}\cdot
  \begin{bmatrix}y_1' &\cdots &y_n'  \end{bmatrix}
\]

With the base case settled, assume $m\geq 3$ and the result holds
for matrices with $m-1$ rows.  As one has $\detrank R(B)\leq 1$, the
induction hypothesis yields matrices $B'$ and $B''$, of size
$(m-1)\times 1$ and $1\times n$ respectively, such that $B=B'B''$.
Set
\[
B'''= \begin{bmatrix}
b'_1 & 0 \\
\vdots & \vdots\\
b'_{m-1} & 0 \\
0 & 1
\end{bmatrix}
\quad \text{and} \quad C=\begin{bmatrix}
b''_1 & \cdots & b''_n \\
a_{m1} & \cdots & a_{mn}
\end{bmatrix}
\]
Evidently, $A= B'''C$.  Now $B\ne 0$ implies $B'\ne 0$, so one gets $\detrank
R(C')\leq 1$.  The basis of the induction shows $C=C'C''$ for a $2\times
1$ matrix $C'$ and a $1\times n$ matrix $C''$. The matrices $A'=B'''C'$
and $A''=C''$ provide the desired decomposition.

Let now $R$ be a domain whose principal ideals satisfy the maximality
condition and over which every $2\times 2$ matrix with zero determinant
has inner rank $1$.  To prove that $R$ is factorial it suffices to
show that for all $a,b\in R$ the ideal $(a)\cap(b)$ is principal, see
\cite[(0.9.4)]{Co-book}.  By hypothesis, there is an $u\in R$ such that
$(ua)$ is maximal among principal ideals in $(a)\cap(b)$.  We are going
to prove $(a)\cap(b)=(ua)$.

Indeed, for each element $xa\in(a)\cap(b)$ there is a unique $x'\in R$
satisfying $xa=x'b$.  In view of the equalities $x/x'=b/a=u/u'$, the
first matrix below has determinantal rank $1$, so the hypothesis on the
ring yields a factorization
 \[
\begin{bmatrix}
u & u' \\
x & x'
\end{bmatrix}
=
\begin{bmatrix}
c \\ d
\end{bmatrix}
\begin{bmatrix}
y & z
\end{bmatrix}
=
\begin{bmatrix}
cy & cz \\
dy & dz
\end{bmatrix}
 \]
Using the first rows of the matrices one gets $cya=ua=u'b=czb$,
hence $ya=zb$.  Thus, there are inclusions of ideals
$(ua)\subseteq(ya)=(zb)\subseteq(a)\cap(b)$; the maximality of
$(ua)$ implies $(ua)=(ya)$.  Using this equality and the second
rows of the matrices, one now obtains $xa=dya\in(dua)\subseteq(ua)$,
as desired.
 \end{proof}

 \section*{Acknowledgements}

We thank Lars Winther Christensen, Claudia Miller, and Greg Piepmeyer
for discussions on this paper and for collaborations
on the related projects \cite{ABCIP}, \cite{ABIM}.



\begin{thebibliography}{23}

\bibitem{AP}
C.~Allday, V.~Puppe,
\textit{Cohomological methods in transformation groups},
Cambridge Stud. Adv. Math. \textbf{32},
Cambridge Univ. Press, Cambridge, 1993.

\bibitem{AB}
L.~L.~Avramov, R.-O.~Buchweitz,
\textit{Lower bounds for Betti numbers},
Compositio Math. \textbf{86} (1993), 147--158.

\bibitem{ABCIP}
L.~L.~Avramov, R.-O.~Buchweitz, L.~W.~Christensen, S.~Iyengar, G.~Piepmeyer,
\textit{Homotopical algebra of differential modules},
in preparation.

\bibitem{ABIM}
L.~L.~Avramov, R.-O.~Buchweitz, S.~Iyengar, C.~Miller,
\textit{Homology of perfect complexes}, 
{\tt http://arxiv.org/abs/math.AC/0609008}.

\bibitem{BE}
D.~Buchsbaum, D.~Eisenbud,
\textit{Algebra structures for finite free resolutions, and
some structure theorems for ideals of codimension $3$},
Amer. J. Math. \textbf{99} (1977), 447--485.

\bibitem{BH}
W.~Bruns and J.~Herzog,
\textit{Cohen-Macaulay rings} (Revised ed.),
Cambridge Stud. Adv. Math. \textbf{39},
Cambridge Univ. Press, Cambridge, 1998.

\bibitem{Ca:inv}
G.~Carlsson,
\textit{On the homology of finite free $(\BZ/2)^k$-complexes},
Invent. Math. \textbf{74} (1983), 139--147.

\bibitem{Ca:ln}
G.~Carlsson,
\textit{Free $(\BZ/2)^k$-actions and a problem in commutative algebra},
Transformation Groups (Pozna\'n, 1985),
Lecture Notes Math. \textbf{1217},
Springer, Berlin, 1986; 79--83.

\bibitem{Ca:ams} G.~Carlsson,
\textit{Free $(\BZ/2)^k$-actions on finite complexes}, Algebraic
Topology and Algebraic K-theory (Princeton, N.J., 1983), Ann. of
Math. Stud. \textbf{113}, Princeton Univ. Press, Princeton, NJ,
1987; 332--344.

\bibitem{CE} H.~Cartan, S.~Eilenberg,
\textit{Homological algebra}, Princeton Univ. Press, Princeton,
NJ, 1956.

\bibitem{Co-book}
P.~M.~Cohn,
\textit{Free ideal rings and localization in general rings},
New Math. Monographs \textbf{3}, Cambridge Univ. Press,
Cambridge, 2006.

\bibitem{DGI}
W.~G.~Dwyer, J.~P.~C.~Greenlees, S.~Iyengar,
\textit{Finiteness in derived categories of local rings},
Commentarii Math. Helvetici, \textbf{81} (2006), 383--432.

\bibitem{Ei}
D.~Eisenbud, \textit{Homological algebra on a complete
intersection, with an application to group representations},
Trans. Amer. Math. Soc. \textbf{260} (1980), 35--64.

\bibitem{EG}
E.~G.~Evans, P.~Griffith,
\textit{The syzygy problem},
Ann. of Math. (2) \textbf{114} (1981), 323--333.

\bibitem{Hl}
S.~Halperin,
\textit{Rational homotopy and torus actions},
Aspects of topology,
London Math. Soc. Lecture Note Ser., \textbf{93},
Cambridge Univ. Press, Cambridge, 1985; 293--306.

\bibitem{Ha}
R.~Hartshorne,
\textit{Algebraic vector bundles on projective spaces:
A problem list},
Topology \textbf{18} (1979), 117--128.

\bibitem{Ho:LNM}
M.~Hochster,
\textit{Cohen-Macaulay modules},
Conference on Commutative Algebra (Lawrence, 1972),
Lecture Notes Math. \textbf{311},
Springer, Berlin, 1973; 120--152.

\bibitem{Ho:CBMS}
M.~Hochster,
\textit{Topics in the homological theory of modules over
commutative rings},
Conf. Board Math. Sci. \textbf{24},
Amer. Math. Soc., Providence, RI, 1975.

\bibitem{Ho:JA} M.~Hochster,
\textit{Big Cohen-Macaulay algebras in dimension three via
Heitmann's theorem},
J. Algebra \textbf{254} (2002), 395--408.

\bibitem{Rb:Sem}
P.~Roberts, \textit{Homological invariants of modules over
commutative rings}, Sem. Math. Sup. \textbf{72}, Presses Univ.
Montr\'eal, Montr\'eal, 1980.

\bibitem{Rb:Camb}
P.~Roberts,
\textit{Multiplicities and Chern classes in local algebra},
Cambridge Tracts Math. \textbf{133},
Cambridge Univ. Press, Cambridge, 1998.
 \end{thebibliography}
 \end{document}